\documentclass[thmsa,onecolumn,10pt,leqno]{article}
\usepackage{amssymb}
\usepackage{amsfonts}
\usepackage{graphicx}
\usepackage{amsmath}

\newtheorem{theorem}{Theorem}

\newtheorem{lemma}[theorem]{Lemma}

\newtheorem{proposition}[theorem]{Proposition}

\newenvironment{proof}[1][Proof]{\textbf{#1.} }{\ \rule{0.5em}{0.5em}}
\newcommand{\R}{\mathbb{R}}
\newcommand{\egl}{\stackrel{\text{(law)}}{\ =\ }}

\begin{document}

\title{Decompositions of stochastic processes based on irreducible group
representations}
\author{Giovanni PECCATI\thanks{%
Laboratoire de Statistique Th\'{e}orique et Appliqu\'{e}e,
Universit\'{e} Paris VI, France. E-mail:
\texttt{giovanni.peccati@gmail.com}} \ \ and \ Jean-Renaud
PYCKE\thanks{Department of mathematics, University of \'{E}vry,
France. E-mail: \texttt{jrpycke@maths.univ-evry.fr}}}
\date{September 23, 2005}
\maketitle

\begin{abstract}
Let $G$ be a topological compact group acting on some space $Y$. We study a
decomposition of $Y$-indexed stochastic processes, based on the
orthogonality relations between the characters of the irreducible
representations of $G$. In the particular case of a Gaussian process with a $%
G$-invariant law, such a decomposition gives a very general
explanation of a classic identity in law -- between quadratic
functionals of a Brownian bridge -- due to Watson (1961). Several
relations with Karhunen-Lo\`{e}ve expansions are discussed, and
some applications and extensions are given -- in particular
related to Gaussian processes indexed by a torus.

\textbf{Key words -- }Stochastic processes; Topological compact
groups; Irreducible representations; Quadratic functionals;
Watson's duplication\ identity; double Wiener-It\^{o} integrals.
\textbf{MSC Classification --} 60G15 60G07 60H05 20C15
\end{abstract}

\section{Introduction}

Let $G$ be a topological compact group acting on a set $Y$, and let $Z\left(
\omega ,y\right) =Z\left( y\right) $ be a stochastic process indexed by the
elements of $Y$. In this paper, we study a decomposition of the paths of $Z$%
, realized by means of the orthogonality relations between the characters of
the irreducible representations of $G$ (see \cite{DuiKolk} or \cite{Serre}
for any unexplained notion about representations). More specifically, we
define $L^{2}\left( G\right) $ to be the space of complex-valued functions
on $G$ that are square integrable w.r.t. the Haar measure, and we note $%
\widehat{G}$ the dual of $G$ (i.e., $\widehat{G}$ is the collection of the
equivalence classes of irreducible representations of $G$). Then, a classic
result of representation theory states that $L^{2}\left( G\right) $ can be
decomposed into an orthogonal sum of finite dimensional spaces, indexed by
the elements $\left[ \pi \right] $ of $\widehat{G}$ and known as the \textit{%
spaces of} \textit{matrix coefficients }of irreducible representations. The
projection operators on such orthogonal spaces have the form of convolutions
with respect to the corresponding characters. Now note $\left( g,y\right)
\mapsto g\cdot y$ the action of $G$ on $Y$, and consider a stochastic
process $Z\left( y\right) $, $y\in Y$, such that, for every fixed $y\in Y$,
the application $g\mapsto Z\left( g\cdot y\right) $ is in $L^{2}\left(
G\right) $. One of our main results states that, if the law of $Z$ is
\textit{invariant} with respect to the action of $G$, the above described
decomposition of $L^{2}\left( G\right) $ translates into a (unique)
decomposition of $Z$ into the sum of simpler stochastic processes, each
indexed by a distinct element of $\widehat{G}$. We write $Z=\sum_{\left[ \pi %
\right] \in \widehat{G}}Z^{\pi }$ for such a decomposition. In Section 3, we
shall prove that, if $Z$ has a $G$-invariant law, then, for distinct $\left[
\pi \right] ,\left[ \sigma \right] \in \widehat{G}$, the processes $Z^{\pi }$
and $Z^{\sigma }$ are non-correlated (in a probabilistic sense), and such
that their paths are orthogonal with respect to \textit{any} $G$-invariant
measure on the parameter space $Y$. In particular, when $Z$ is Gaussian and $%
\left[ \pi \right] $ and $\left[ \sigma \right] $ have real characters, $%
Z^{\pi }$ and $Z^{\sigma }$ are also Gaussian, and therefore stochastically
independent. In the last section we shall discuss some connections between
our decomposition and the Karhunen-Lo\`{e}ve expansion (see for instance
\cite{Adler}) of suitably regular Gaussian processes.

\bigskip

As discussed below, the initial impetus for such an investigation was
provided by the following \textit{duplication identity} due to Watson (see
\cite{Wat}, and \cite{PeYor2005} for two-parameter generalizations): if $b$
is a standard Brownian bridge on $\left[ 0,1\right] $, from 0 to 0, then

\begin{equation}
\int_{0}^{1}\left( b\left( s\right) -\int_{0}^{1}b\left( u\right) du\right)
^{2}ds\overset{law}{=}\frac{1}{4}\left\{ \int_{0}^{1}b\left( s\right)
^{2}ds+\int_{0}^{1}b_{\ast }\left( s\right) ^{2}ds\right\} ,
\label{Watson1961}
\end{equation}%
where $b_{\ast }$ is an independent copy of $b$. The reader is referred to
\cite[p. 220]{SW} for a proof of (\ref{Watson1961}) using Karhunen-Lo\`{e}ve
expansions, and to \cite{ShiYor} for a probabilistic discussion based on
several identities in law between Brownian functionals. However, the short
proof of (\ref{Watson1961}) recently given by the second author (see \cite%
{Pycke}) suggests that there is a simple algebraic structure behind such a
duplication result. As a by-product of our analysis, we will indeed show
that (\ref{Watson1961}) derives from a very particular case of the
decomposition described above. In particular, our results will make clear
that there are two crucial elements behind (\ref{Watson1961}), namely: (i)
since $b_{t}\overset{law}{=}b_{1-t}$ (as stochastic processes), the law of $%
b\left( \cdot \right) -\int_{0}^{1}b\left( u\right) du$ is invariant with
respect to the elementary action, of $G=\left\{ 1,g\right\} \simeq \mathbb{Z}%
/2\mathbb{Z}$ on $\left[ 0,1\right] $, given by $1\cdot t=t$ and $g\cdot
t=1-t$, and (ii) Lebesgue measure is invariant with respect to the same
action of $G$. It follows that the above described theory can be applied,
and (\ref{Watson1961}) turns out to be the result of an orthogonal
decomposition of the paths of $b\left( \cdot \right) -\int_{0}^{1}b\left(
u\right) du$ into two independent components. More to the point, the
(somewhat mysterious) factor $1/4$ on the right hand side of (\ref%
{Watson1961}) appears as the square of a \textit{normalization factor }($%
1/4=1/\left\vert G\right\vert ^{2}$), which enters quite naturally into the
expression of the projection operators associated to matrix coefficients.
The generalizations of (\ref{Watson1961}) given in \cite{PeYor2005} have
similar interpretations in terms of group representations (see Section 3
below for the analysis of a \textit{quadruplication identity}).

\bigskip

Although we are mainly motivated by finite groups and Gaussian processes, we
will work in a very general framework, allowing $G$ to be compact, and
considering processes that are not necessarily Gaussian (see for instance
Paragraph 3.1).

\bigskip

The paper is organized as follows. In Section 2 we recall some basic facts
about group representations and related orthogonality relations. Section 3
deals with decompositions of stochastic processes, first in the general
case, and then (Section 3.2) in the specific setup of Gaussian processes. In
Section 4, we study generalized (Gaussian) Volterra processes, and establish
some necessary and sufficient conditions (based on the method of cumulants)
to have that such processes verify a relation analogous to (\ref{Watson1961}%
). In the last Section, we discuss several refinements and
applications, mainly related to Karhunen-Lo\`{e}ve expansions,
Gaussian processes indexed by a torus.

\section{Preliminaries and main results from group representation theory}

\subsection{Representations of compact groups and orthogonal decompositions}

In this section, we shall present several definitions and results from the
theory of representations of topological compact groups. Our use of this
theory is mainly inspired by the discussion contained in \cite[Chapter IV]%
{DuiKolk}, where a strong accent is placed on the so-called \textbf{%
Peter-Weyl theorem }(see \cite[Theorem 4.6.1]{DuiKolk}, as well as the
discussion below), and its consequences in terms of the decomposition of the
$L^{2}$ space associated to a topological compact group (when endowed with
its Haar measure). The reader is referred to \cite{DuiKolk} for any
unexplained definition or result. Other references for group representations
are the classic monographs \cite{Serre} and \cite{Diaconis}.

\bigskip

A \textbf{topological group }is a pair $\left( G,\mathbb{G}\right) $, where $%
G$ is a group and $\mathbb{G}$ is a topology such that the following three
conditions are satisfied: (i) $G$ is a Hausdorff topological space, (ii) the
multiplication $G\times G\mapsto G:\left( g,h\right) \mapsto gh$ is
continuous, (iii) the inversion $G\mapsto G:g\mapsto g^{-1}$ is continuous.
In what follows, when no further specification is given, $G$ will always
denote a topological group (the topology $\mathbb{G}$ being implicitly
defined) which is also \textbf{compact }(see e.g. \cite[p. 34]{Dudley}) and
such that $\mathbb{G}$ has a countable basis. For such a $G$, we will denote
by $C\left( G\right) $ the class of continuous, complex-valued functions on $%
G$; $\mathcal{G}$ is the (Borel) $\sigma $-field generated by $\mathbb{G}$.
An immediate consequence (see \cite[Section 10.3]{DuiKolk}) of the structure
imposed on $G$, is that $G$ always carries a (unique) positive Borel
measure, noted $dg$ and known as the \textbf{Haar measure}, such that $%
\int_{G}dg=1$, and $\forall f\in C\left( G\right) $ and $\forall h\in G$%
\begin{eqnarray*}
\int_{G}f\left( g\right) dg &=&\int_{G}f\left( g^{-1}\right) dg \\
\int_{G}f\left( hg\right) dg &=&\int_{G}f\left( gh\right) dg=\int_{G}f\left(
g\right) dg\text{ \ \ (left and right invariance);}
\end{eqnarray*}%
we shall note $L^{2}\left( G,dg\right) =L^{2}\left( G\right) $ the Hilbert
space of complex valued functions on $G$ that are square integrable with
respect to $dg$, endowed with the usual inner product $\left\langle
f_{1},f_{2}\right\rangle _{G}=\int_{G}f_{1}\left( g\right) \overline{%
f_{2}\left( g\right) }dg$. We note $\left\Vert \cdot \right\Vert _{G}$ the
norm associated to $\left\langle \cdot ,\cdot \right\rangle _{G}$, and we
observe that $L^{2}\left( G\right) $ is the completion of $C\left( G\right) $
with respect to $\left\Vert \cdot \right\Vert _{G}.$

\bigskip

\textbf{Remark -- }When $G$ is finite, then $\mathbb{G}$ is necessarily the
discrete topology, and $dg$ coincides with the \textbf{normalized counting
measure} associated to $G$, that is%
\begin{equation*}
dg=\frac{1}{\left\vert G\right\vert }\sum_{h\in G}\delta _{h}\left(
dg\right) \text{,}
\end{equation*}%
where $\delta _{h}\left( \cdot \right) $ stands for the Dirac mass
concentrated at $h$, and $\left\vert G\right\vert $ is the cardinality of $G$%
.

\bigskip

Let $V$ be a topological vector space over $\mathbb{C}$. A \textbf{%
representation }of $G$ in $V$ is an homomorphism $\pi $, from $G$ into $%
\mathbf{GL}\left( V\right) $ (the set of complex isomorphisms of $V$ into
itself), such that the mapping $G\times V\mapsto V:\left( g,v\right) \mapsto
\pi \left( g\right) \left( v\right) $ is continuous. The \textbf{dimension }$%
d_{\pi }$ of a representation $\pi $ is defined to be the dimension of $V$.
A representation $\pi $ of $G$ in $V$\textbf{\ }is \textbf{irreducible}, if
the only closed $\pi \left( G\right) $-invariant subspaces of $V$ are $%
\left\{ 0\right\} $ and $V$. It is well known that irreducible
representations are defined up to equivalence classes (see \cite[p. 210]%
{DuiKolk}). Following \cite{DuiKolk}, we will denote by $\left[ \pi \right] $
the equivalence class of a given irreducible representation $\pi $; the set
of equivalence classes of irreducible representations of $G$ is noted $%
\widehat{G}$, and it is called the \textbf{dual }of $G$. Note that, in our
setting, irreducible representations\textit{\ are always finite dimensional}%
. Moreover, we will systematically assume (without loss of generality, see
\cite[Corollary 4.2.2]{DuiKolk}) that every irreducible representation is
also unitary. Finally, we recall that, according e.g. to \cite[Theorem 4.3.4
(v)]{DuiKolk}, since $G$ is second countable $\widehat{G}$ is necessarily
countable.

\bigskip

To every finite dimensional representation $\pi :G\mapsto \mathbf{GL}\left(
V\right) $ we associate the mapping%
\begin{equation*}
\chi _{\pi }:G\mapsto \mathbb{C}:g\mapsto \text{Trace }\pi \left( g\right)
\text{,}
\end{equation*}%
called the \textbf{character }of $\pi $. Two finite dimensional
representations are equivalent if, and only if, they have the same
character. Moreover, it is easily seen that characters are \textit{central%
\footnote{%
That is, for every $x,g\in G$, $\chi _{\pi }\left( x^{-1}gx\right) =\chi
_{\pi }\left( g\right) $.}}, and \textit{continuous} functions on $G$.

\bigskip

In this paper, we shall develop some Hilbert space techniques that are
directly based on the orthogonality relations between the characters of
distinct irreducible representations. To this end, for every $\left[ \pi %
\right] \in \widehat{G}$ we associate a finite dimensional subspace $M_{\pi
}\subseteq L^{2}\left( G\right) $ in the following way. Select an element $%
\pi :G\mapsto \mathbf{GL}\left( V\right) $ in $\left[ \pi \right] $, as well
as a basis $\mathbf{e}=\left\{ e_{1},...,e_{n}\right\} $ of $V$ (plainly, $%
n=d_{\pi }$) with respect to which $\pi $ is unitary; the space $M_{\pi }$
is defined as the set of the (complex) linear combinations of the \textbf{%
matrix coefficients} associated to $\pi $ and to the basis $\mathbf{e}$,
that is, $M_{\pi }$ is composed by the linear combinations of the functions%
\begin{equation*}
g\mapsto \pi \left( g\right) _{k}^{j}\text{, \ \ }j,k=1,...,n\text{,}
\end{equation*}%
where, for each $g\in G$, $\left\{ \pi \left( g\right)
_{k}^{j}:j,k=1,...,n\right\} $ is the matrix representation of $\pi \left(
g\right) $ with respect to the basis $\mathbf{e}$. Note that such a
definition is well given, as $M_{\pi }$ does not depend on the choice of the
representative element of $\left[ \pi \right] $ and of the basis of $V$. Of
course, $M_{\pi }$ is finite dimensional (and therefore closed; more
precisely: $\dim M_{\pi }=d_{\pi }^{2}$, see \cite[Theorem 4.3.4]{DuiKolk})
and $M_{\pi }\subseteq C\left( G\right) $, for every $\left[ \pi \right] \in
\widehat{G}$.

\bigskip

Before stating one of the crucial results for our analysis, we introduce a
\textbf{convolution }operation on $L^{2}\left( G\right) $, which is defined,
for $f,k\in L^{2}\left( G\right) $, by the formula%
\begin{equation}
\left( f\ast k\right) \left( u\right) =\int_{G}f\left( g\right) k\left(
g^{-1}u\right) dg=\int_{G}f\left( ug^{-1}\right) k\left( g\right) dg\text{,
\ \ }u\in G.  \label{conv}
\end{equation}

The following result summarizes all the orthogonality relations --
associated to the notion of character -- that are relevant to our discussion
(for proofs and further analysis in this direction, the reader is referred
to \cite[paragraphs 4.2--4.6]{DuiKolk}).

\begin{theorem}
Let the above notation and assumptions prevail. Then,

\begin{enumerate}
\item if $\left[ \pi \right] ,\left[ \sigma \right] \in \widehat{G}$, and $%
\left[ \pi \right] \neq \left[ \sigma \right] $, the spaces $M_{\pi }$ and $%
M_{\sigma }$ are orthogonal in $L^{2}\left( G\right) ;$

\item for every $\left[ \pi \right] \in \widehat{G}$ the orthogonal
projection operator, from $L^{2}\left( G\right) $ to $M_{\pi }$, is given by%
\begin{equation}
E_{\pi }:L^{2}\left( G\right) \mapsto M_{\pi }:f\mapsto d_{\pi }\left( f\ast
\chi _{\pi }\right) :=E_{\pi }f;  \label{proj}
\end{equation}

\item the class $\left\{ M_{\pi }:\left[ \pi \right] \in \widehat{G}\right\}
$ is total in $L^{2}\left( G\right) $, and therefore%
\begin{equation}
L^{2}\left( G\right) =\bigoplus\limits_{\left[ \pi \right] \in \widehat{G}%
}M_{\pi }  \label{generalortho}
\end{equation}%
where $\bigoplus $ stands for a direct Hilbert space (orthogonal) sum;

\item for every $\left[ \pi \right] \in \widehat{G}$,
\begin{equation}
\left\langle \chi _{\pi },\chi _{\pi }\right\rangle _{G}=1\text{, \ }E_{\pi
}\chi _{\pi }=\chi _{\pi }\text{ \ \ and \ }E_{\sigma }\chi _{\pi }=0\text{
\ \ if }\left[ \sigma \right] \in \widehat{G}\text{\ and\ }\left[ \pi \right]
\neq \left[ \sigma \right] ,  \label{charactersortho}
\end{equation}%
and consequently $\left\{ \chi _{\pi }:\left[ \pi \right] \in \widehat{G}%
\right\} $ is an orthonormal system in $L^{2}\left( G\right) $.
\end{enumerate}
\end{theorem}

\textbf{Remarks -- }(i) Since, for every $\left[ \pi \right] \in \widehat{G}$%
, the function $G\ni g\mapsto \chi _{\pi }\left( g\right) $ is
conjugacy-invariant, the projection $E_{\pi }f$, as defined in (\ref{proj}),
is also equal to $d_{\pi }\left( \chi _{\pi }\ast f\right) .$

(ii) Point 3 of Theorem 1 can be seen as a direct consequence of the \textbf{%
Peter-Weyl theorem} (see \cite[Theorem 4.6.1]{DuiKolk}), stating that the
space of the linear combinations of the matrix coefficients, associated to
the finite-dimensional representations of $G$, is dense in $C\left( G\right)
$, endowed with the supremum norm.

(iii) For future reference, we recall that the following four conditions are
equivalent (see \cite[p. 235]{DuiKolk}): (a) $G$ is Abelian, (b) $d_{\pi }=1$
for every $\left[ \pi \right] \in \widehat{G}$, (c) every $f\in C\left(
G\right) $ is conjugacy-invariant, (d) the convolution operation defined in (%
\ref{conv}) is commutative. In particular, if $G$ is Abelian, then the
system $\left\{ \chi _{\pi }:\left[ \pi \right] \in \widehat{G}\right\} $ is
orthonormal \textit{and complete} in $L^{2}\left( G\right) $. If $G$ is
Abelian and finite, then $\left\vert G\right\vert =$ $\mid \widehat{G}\mid $.

\subsection{Actions and decompositions of complex-valued functions}

Consider a measurable space $\left( X,\mathcal{X}\right) $. In this paper, a
\textbf{left} \textbf{action }$A$ of $G$ on $X$ is a $\mathcal{G}\otimes
\mathcal{X}$ -- measurable function, from $G\times X$ to $X$ (recall that $%
\mathcal{G}$ is the Borel $\sigma $-field of $G$), such that, for every $%
g,h\in G$ and $x\in X$,%
\begin{equation*}
A\left( gh,x\right) =A\left( g,A\left( h,x\right) \right) .
\end{equation*}

A right action can be defined in a similar way, but we will deal only with
left actions; for the sake of simplicity, in the sequel left actions are
simply called \textbf{actions}. When there is no ambiguity on the action $A$%
, we will sometimes use the customary abbreviation%
\begin{equation*}
A\left( g,x\right) =g\cdot x\text{ \ \ (}g\in G\text{, }x\in X\text{).}
\end{equation*}

A $\sigma $-finite, positive measure $\nu $ on $\left( X,\mathcal{X}\right) $%
, is said to be \textbf{invariant with respect to the action} $A$ \textbf{of}
$G$ (or simply, again when there is no ambiguity on the action $A$, $G$-%
\textbf{invariant}) if, for every complex valued function $f\in L^{1}\left(
\nu \right) $,
\begin{equation*}
\int_{X}f\left( x\right) \nu \left( dx\right) =\int_{X}f\left( A\left(
g,x\right) \right) \nu \left( dx\right) =\int_{X}f\left( g\cdot x\right) \nu
\left( dx\right) ,
\end{equation*}%
for every $g\in G$.

\bigskip

Now fix an action $A$ of $G$ on $X$, and consider a measurable function $%
Z:X\mapsto \mathbb{C}$. We associate to $Z$ the function
\begin{equation}
Z_{\ast }:G\times X\mapsto \mathbb{C}:\left( g,x\right) \mapsto Z_{\ast
}\left( g,x\right) =Z\left( A\left( g,x\right) \right) =Z\left( g\cdot
x\right) ,  \label{jointfunct}
\end{equation}%
which is of course $\mathcal{G}\otimes \mathcal{X}$ -- measurable. For each
fixed $x\in X$, we define the $\mathcal{G}$ -- measurable function
\begin{equation}
Z_{G}\left[ x\right] :G\mapsto \mathbb{C}:g\mapsto Z_{\ast }\left(
g,x\right) \text{ ;}  \label{marg1}
\end{equation}%
analogously, for each fixed\textit{\ }$g\in G$, we note%
\begin{equation}
Z_{X}\left[ g\right] :X\mapsto \mathbb{C}:x\mapsto Z_{\ast }\left(
g,x\right) \text{,}  \label{marg2}
\end{equation}%
which defines in turn a $\mathcal{X}$ -- measurable mapping. If, for some
fixed $x\in X$, the above introduced function $Z_{G}\left[ x\right] $ is an
element of $L^{2}\left( G\right) $, we set, for each $\left[ \pi \right] \in
\widehat{G}$ and $g\in G$,
\begin{equation}
Z_{\ast }^{\pi }\left( g,x\right) =E_{\pi }Z_{G}\left[ x\right] \left(
g\right) ,  \label{projmarg1}
\end{equation}%
where, by using Theorem 1 and (\ref{conv}),
\begin{equation}
E_{\pi }Z_{G}\left[ x\right] \left( g\right) =d_{\pi }\int_{G}\chi _{\pi
}\left( h\right) Z_{G}\left[ x\right] \left( h^{-1}g\right) dh=d_{\pi
}\int_{G}\chi _{\pi }\left( h\right) Z\left( h^{-1}g\cdot x\right) dh,\text{
\ \ }g\in G\text{.}  \label{projoperator}
\end{equation}

As usual, we write $e$ to indicate the \textbf{identity }element of the
group $G$. If, for a measurable mapping $Z:X\mapsto \mathbb{C}$, $Z_{G}\left[
x\right] \in L^{2}\left( G\right) $ for every $x\in X$, we introduce the two
functions, defined respectively for a fixed $x\in X$ and for a fixed $g\in G$%
,%
\begin{eqnarray}
Z_{G}^{\pi }\left[ x\right] &:&G\mapsto \mathbb{C}:g\mapsto Z_{\ast }^{\pi
}\left( g,x\right)  \label{margmarg} \\
Z_{X}^{\pi }\left[ g\right] &:&X\mapsto \mathbb{C}:x\mapsto Z_{\ast }^{\pi
}\left( g,x\right) ;  \notag
\end{eqnarray}%
finally, for such a $Z$, we set%
\begin{equation}
Z^{\pi }\left( x\right) =Z_{X}^{\pi }\left[ e\right] \left( x\right)
=Z_{\ast }^{\pi }\left( e,x\right) =E_{\pi }Z_{G}\left[ x\right] \left(
e\right) =d_{\pi }\int_{G}\chi _{\pi }\left( g\right) Z\left( g^{-1}\cdot
x\right) dg\text{, \ \ }x\in X.  \label{projmarg2}
\end{equation}

Note that, since $A$ is a left action,
\begin{equation}
Z_{X}^{\pi }\left[ g\right] \left( x\right) =Z^{\pi }\left( g\cdot x\right) ,%
\text{ \ \ }g\in G\text{, \ }x\in X.  \label{projmarg3}
\end{equation}

\bigskip

\textbf{Remark -- }If the function\textbf{\ }$Z$ is such that $Z_{G}\left[ x%
\right] \in L^{2}\left( G\right) $ for every $x\in X$, then the mapping $%
\left( g,x\right) \mapsto Z_{\ast }^{\pi }\left( g,x\right) $ is $\mathcal{G}%
\otimes \mathcal{X}$ -- measurable. It follows that the two mappings $%
Z_{G}^{\pi }\left[ x\right] $ and $Z_{X}^{\pi }\left[ g\right] $ defined in (%
\ref{margmarg}) are, respectively, $\mathcal{G}$ -- measurable and $\mathcal{%
X}$ -- measurable. In particular, the application $x\mapsto Z^{\pi }\left(
x\right) $ (as defined in (\ref{projmarg2})) is a $\mathcal{X}$ --
measurable mapping.

\bigskip

The following result turns out to be the key tool of our analysis.

\begin{theorem}
Under the above notation and assumptions, fix an action $A$ of the group $G$
on $X$. Consider moreover two measurable functions $S,Z:X\mapsto \mathbb{C}$%
, such that for each $x\in X$, $S_{G}\left[ x\right] ,Z_{G}\left[ x\right]
\in L^{2}\left( G\right) $. Then,

\begin{enumerate}
\item for any $\left[ \pi \right] ,\left[ \sigma \right] \in \widehat{G}$
such that $\left[ \pi \right] \neq \left[ \sigma \right] $, and for
arbitrary points $x_{1},x_{2}\in X$, the following orthogonality relation is
satisfied:%
\begin{equation}
\left\langle S_{G}^{\pi }\left[ x_{1}\right] ,Z_{G}^{\sigma }\left[ x_{2}%
\right] \right\rangle _{G}=0;  \label{ortho1}
\end{equation}

\item for every $x\in X$,
\begin{equation}
Z_{G}\left[ x\right] =\sum_{\left[ \pi \right] \in \widehat{G}}Z_{G}^{\pi }%
\left[ x\right] \text{ \ \ and \ \ }S_{G}\left[ x\right] =\sum_{\left[ \pi %
\right] \in \widehat{G}}S_{G}^{\pi }\left[ x\right] ,  \label{decomposition1}
\end{equation}%
where the convergence of the (possibly infinite) series takes place in $%
L^{2}\left( G\right) ,$ and for any $x_{1},x_{2}\in X$%
\begin{equation}
\left\langle S_{G}\left[ x_{1}\right] ,Z_{G}\left[ x_{2}\right]
\right\rangle _{G}=\sum_{\left[ \pi \right] \in \widehat{G}}\left\langle
S_{G}^{\pi }\left[ x_{1}\right] ,Z_{G}^{\pi }\left[ x_{2}\right]
\right\rangle _{G},  \label{fouriercoeff1}
\end{equation}%
with convergence in $\ell ^{2}$;

\item in addition to the previous assumptions, suppose there exists a $G$%
-invariant measure $\nu $ on $\left( X,\mathcal{X}\right) $, such that the
functions $Z_{\ast }$ and $S_{\ast }$, defined according to (\ref{jointfunct}%
), are elements of
\begin{equation*}
L^{2}\left( G\times X,\mathcal{G}\otimes \mathcal{X},dg\times \nu \left(
dx\right) \right) :=L^{2}\left( dg\times \nu \left( dx\right) \right)
\end{equation*}%
and also, for every $g\in G$, $Z_{X}\left[ g\right] ,S_{X}\left[ g\right]
\in L^{2}\left( X,\mathcal{X},\nu \left( dx\right) \right) $ $:=L^{2}\left(
\nu \left( dx\right) \right) $; then, for every $\left[ \pi \right] \in
\widehat{G}$ and every $g\in G$, $Z_{\ast }^{\pi },S_{\ast }^{\pi }\in
L^{2}\left( dg\times \nu \left( dx\right) \right) $, $Z_{X}^{\pi }\left[ g%
\right] ,S_{X}^{\pi }\left[ g\right] \in L^{2}\left( \nu \left( dx\right)
\right) $, and moreover
\begin{equation}
\int_{X}S_{X}^{\pi }\left[ g\right] \left( x\right) \overline{Z_{X}^{\sigma }%
\left[ g\right] \left( x\right) }\nu \left( dx\right) =0  \label{ortho2}
\end{equation}%
for every $\left[ \pi \right] ,\left[ \sigma \right] \in \widehat{G},$ such
that $\left[ \pi \right] \neq \left[ \sigma \right] $;

\item under the assumptions and notation of point 3.,
\begin{equation}
Z_{\ast }=\sum_{\left[ \pi \right] \in \widehat{G}}Z_{\ast }^{\pi }\text{ \
\ and \ \ }S_{\ast }=\sum_{\left[ \pi \right] \in \widehat{G}}S_{\ast }^{\pi
}  \label{twovardecomp}
\end{equation}%
where the series are orthogonal and convergent in $L^{2}\left( dg\times \nu
\left( dx\right) \right) $, and therefore
\begin{equation}
\left\langle Z_{\ast },S_{\ast }\right\rangle _{L^{2}\left( dg\times \nu
\left( dx\right) \right) }=\sum_{\left[ \pi \right] \in \widehat{G}%
}\left\langle Z_{\ast }^{\pi },S_{\ast }^{\pi }\right\rangle _{L^{2}\left(
dg\times \nu \left( dx\right) \right) };  \label{fourcoeff2}
\end{equation}

\item under the assumptions and notation of point 3., for every $g\in G$,
\begin{equation}
Z_{X}\left[ g\right] \left( x\right) =\sum_{\left[ \pi \right] \in \widehat{G%
}}Z_{X}^{\pi }\left[ g\right] \left( x\right) \text{ \ \ and \ \ }S_{X}\left[
g\right] \left( x\right) =\sum_{\left[ \pi \right] \in \widehat{G}%
}S_{X}^{\pi }\left[ g\right] \left( x\right)  \label{ortho3}
\end{equation}%
where the series are orthogonal and convergent in $L^{2}\left( \nu \left(
dx\right) \right) $, and
\begin{equation}
\left\langle Z_{X}\left[ g\right] ,S_{X}\left[ g\right] \right\rangle
_{L^{2}\left( \nu \left( dx\right) \right) }=\sum_{\left[ \pi \right] \in
\widehat{G}}\left\langle Z_{X}^{\pi }\left[ g\right] ,S_{X}^{\pi }\left[ g%
\right] \right\rangle _{L^{2}\left( \nu \left( dx\right) \right) }.
\label{fourcoeff3}
\end{equation}
\end{enumerate}
\end{theorem}

\bigskip

\textbf{Remark --} When $G$ is finite, $\widehat{G}$ is also finite (since $%
\mid \widehat{G}\mid $ coincides with the number of conjugacy classes in $G$%
). In this case, Theorem 2-2 gives a decomposition of the function $%
Z:X\mapsto \mathbb{C}$. As a matter of fact, for every $x\in X$,%
\begin{equation}
Z\left( x\right) =\sum_{\left[ \pi \right] \in \widehat{G}}Z^{\pi }\left(
x\right) ,  \label{path-decomposition}
\end{equation}%
where the sum is finite, and on the right-hand side we use the notation
introduced in (\ref{projmarg2}).

\bigskip

\textbf{Proof of Theorem 2. }(\textbf{1.}) By definition, $S_{G}^{\pi }\left[
x_{1}\right] $ and $Z_{G}^{\sigma }\left[ x_{2}\right] $ equal the
orthogonal projections, respectively of $S_{G}\left[ x_{1}\right] $ and $%
Z_{G}\left[ x_{2}\right] $, on the finite dimensional spaces $M_{\pi }$ and $%
M_{\sigma }$. Since, according to Theorem 1-1, $M_{\pi }$ and $M_{\sigma }$
are orthogonal in $L^{2}\left( G\right) $, relation (\ref{ortho1}) follows.

\noindent%
(\textbf{2.}) Relation (\ref{decomposition1}) is an immediate consequence of
(\ref{generalortho}), whereas (\ref{fouriercoeff1}) is a standard formula of
the Parseval-Plancherel type.

\noindent%
(\textbf{3.}) Observe first that, by assumption,%
\begin{equation*}
\int_{X}\left[ \int_{G}\left\vert Z_{\ast }\left( h,x\right) \right\vert
^{2}dh\right] \nu \left( dx\right) <+\infty .
\end{equation*}%
Moreover, we observe that, for every $\left[ \pi \right] \in \widehat{G}$,
the continuous function $\left\vert \chi _{\pi }\right\vert :G\mapsto
\mathbb{R}
_{+}$ is bounded by a constant $\alpha _{\pi }\in \left( 0,+\infty \right) $
(since $G$ is compact), and therefore, thanks to the right invariance of the
Haar measure $dh$ and Jensen inequality,
\begin{eqnarray*}
\int_{X}\left\vert Z_{X}^{\pi }\left[ g\right] \left( x\right) \right\vert
^{2}\nu \left( dx\right) &=&d_{\pi }^{2}\int_{X}\left\vert \int_{G}\chi
_{\pi }\left( h\right) Z_{\ast }\left( h^{-1}g,x\right) dh\right\vert
^{2}\nu \left( dx\right) \\
&\leq &d_{\pi }^{2}\alpha _{\pi }^{2}\int_{X}\left[ \int_{G}\left\vert
Z_{\ast }\left( h^{-1}g,x\right) \right\vert ^{2}dh\right] \nu \left(
dx\right) \\
&=&d_{\pi }^{2}\alpha _{\pi }^{2}\int_{X}\left[ \int_{G}\left\vert Z_{\ast
}\left( h,x\right) \right\vert ^{2}dh\right] \nu \left( dx\right) <+\infty .
\end{eqnarray*}%
Also,%
\begin{eqnarray*}
\int_{G}\int_{X}\left\vert Z_{\ast }^{\pi }\left( g,x\right) \right\vert
^{2}\nu \left( dx\right) dg &=&\int_{G}\int_{X}\left\vert Z_{G}^{\pi }\left[
g\right] \left( x\right) \right\vert ^{2}\nu \left( dx\right) dg \\
&\leq &d_{\pi }^{2}\alpha _{\pi }^{2}\int_{G}\int_{X}\left[
\int_{G}\left\vert Z_{\ast }\left( h,x\right) \right\vert ^{2}dh\right] \nu
\left( dx\right) dg \\
&=&d_{\pi }^{2}\alpha _{\pi }^{2}\int_{X}\left[ \int_{G}\left\vert Z_{\ast
}\left( h,x\right) \right\vert ^{2}dh\right] \nu \left( dx\right) <+\infty ,
\end{eqnarray*}%
since $\int_{G}dg=1$. It follows that $Z_{\ast }^{\pi }\in L^{2}\left(
dg\times \nu \left( dx\right) \right) $ and $Z_{X}^{\pi }\left[ g\right] \in
L^{2}\left( \nu \left( dx\right) \right) $ for any $g\in G$, and an
analogous conclusion holds for $S$. We may prove (\ref{ortho2}) by using an
easy declination of the \textquotedblleft averaging\textquotedblright\
technique (see for instance \cite[paragraph 4.2]{DuiKolk}). Indeed, for $%
\left[ \pi \right] ,\left[ \sigma \right] \in \widehat{G}$ such that $\left[
\pi \right] \neq \left[ \sigma \right] $, thanks to formula (\ref{projmarg3}%
) and the $G$-invariance of $\nu $,%
\begin{eqnarray*}
\int_{X}S_{X}^{\pi }\left[ g\right] \left( x\right) \overline{Z_{X}^{\sigma }%
\left[ g\right] \left( x\right) }\nu \left( dx\right) &=&\int_{X}S^{\pi
}\left( g\cdot x\right) \overline{Z^{\sigma }\left( g\cdot x\right) }\nu
\left( dx\right) \\
&=&\int_{G}\int_{X}S^{\pi }\left( h\cdot x\right) \overline{Z^{\sigma
}\left( h\cdot x\right) }\nu \left( dx\right) dh \\
&=&\int_{X}\left[ \int_{G}S^{\pi }\left( h\cdot x\right) \overline{Z^{\sigma
}\left( h\cdot x\right) }dh\right] \nu \left( dx\right) \\
&=&\int_{X}\left\langle S_{G}^{\pi }\left[ x\right] ,Z_{G}^{\sigma }\left[ x%
\right] \right\rangle _{G}\nu \left( dx\right) =0,
\end{eqnarray*}%
where we have used a standard Fubini theorem, as well as Theorem 2-1.

\noindent%
(\textbf{4.}) The first part derives immediately from points 1. and 2., as
well as the fact that $Z_{\ast },S_{\ast }\in L^{2}\left( dg\times \nu
\left( dx\right) \right) $ by assumption. Formula (\ref{fourcoeff2}) is
again of the Parseval-Plancherel type.

\noindent%
(\textbf{5.}) Formula\ (\ref{ortho3}) derives from the elementary relation%
\begin{equation*}
\int_{X}\left\vert Z_{X}\left[ g\right] \left( x\right) -\sum_{\left[ \pi %
\right] \in \widehat{G}}Z_{X}^{\pi }\left[ g\right] \left( x\right)
\right\vert ^{2}\nu \left( dx\right) =\int_{G}\int_{X}\left\vert Z_{\ast
}\left( h,x\right) -\sum_{\left[ \pi \right] \in \widehat{G}}Z_{\ast }^{\pi
}\left( h,x\right) \right\vert ^{2}\nu \left( dx\right) dh=0,
\end{equation*}%
where the first equality is due to the $G$-invariance of $\nu $, and the
second comes from point 4. Relation (\ref{fourcoeff3}) is straightforward. $%
\blacksquare $

\section{Decompositions of stochastic processes}

\subsection{General results}

Let $\left( \Omega ,\mathcal{F},\mathbf{P}\right) $ be a probability space,
and let $\left( Y,\mathcal{Y}\right) $ be a measurable space. A $Y$-\textbf{%
indexed} \textbf{stochastic process }$Z$ is a $\mathcal{F}\otimes \mathcal{Y}
$ -- measurable application $Z:\Omega \times Y\mapsto \mathbb{C}:\left(
\omega ,y\right) \mapsto Z\left( \omega ,y\right) $ \footnote{%
We write $Z\left( y\right) $, for fixed $y\in Y$, to indicate the random
variable $\omega \mapsto Z\left( \omega ,y\right) $.}. To simplify some
arguments, we shall systematically suppose that the $\sigma $-field $%
\mathcal{F}$ contains singletons, that is, $\left\{ \omega \right\} \in
\mathcal{F}$ for every $\omega \in \Omega $.

\bigskip

In this section, the product space $\Omega \times Y$ will play roughly the
same role as the space $\left( X,\mathcal{X}\right) $ in Section 2. As a
consequence, we shall sometimes use the compact notation

\begin{equation}
\Omega \times Y=X_{0}\text{, \ \ }\mathcal{F}\otimes \mathcal{Y}=\mathcal{X}%
_{0},  \label{omegaY}
\end{equation}%
and write $x_{0}$ to indicate the generic element $\left( \omega ,y\right) $
of $X_{0}$. Given a topological compact group $G$ and an action $A$ of $G$
on $X_{0}$, for fixed $y\in Y$ and $g\in G$, we write $Z\left( g\cdot
y\right) $ to indicate the random variable%
\begin{equation*}
\Omega \ni \omega \mapsto Z\left( g\cdot \left( \omega ,y\right) \right) .
\end{equation*}

We say that the law of a family $\mathbf{Z}=\left\{ Z_{i}:i\in I\right\} $
of stochastic processes is \textbf{invariant with respect to the action} $A$
of $G$ on $X_{0}$ (or, simply, $G$-\textbf{invariant}) if, for every $n\geq
1 $ and every continuous, bounded function $f$ on $\mathbb{C}^{n}$%
\begin{equation*}
\mathbf{E}\left[ f\left( Z_{i_{1}}\left( y_{1}\right) ,...,Z_{i_{n}}\left(
y_{n}\right) \right) \right] =\mathbf{E}\left[ f\left( Z_{i_{1}}\left(
g\cdot y_{1}\right) ,...,Z_{i_{n}}\left( g\cdot y_{n}\right) \right) \right]
\end{equation*}%
for every $g\in G$, every $\left( y_{1},...,y_{n}\right) \in Y^{n}$, and
every $\left( i_{1},...,i_{n}\right) \in I^{n}$.

\bigskip

\textbf{Remark -- }Every action $A^{\prime }$ of $G$ on $Y$ always defines
an action $A$ on $X_{0}$, through the relation: for every $x_{0}=\left(
\omega ,y\right) \in X_{0}$,%
\begin{equation}
A\left( g,x_{0}\right) =g\cdot x_{0}=\left( \omega ,A^{\prime }\left(
g,y\right) \right) .  \label{canaction}
\end{equation}%
Analogously, every action $\underline{A}^{\prime }$ of $G$ on $\Omega $
defines an action $\underline{A}$ on $X_{0}$: for every $x_{0}=\left( \omega
,y\right) \in X_{0}$,%
\begin{equation}
\underline{A}\left( g,x_{0}\right) =g\cdot x_{0}=\left( \underline{A}%
^{\prime }\left( g,\omega \right) ,y\right) .  \label{canaction2}
\end{equation}

\bigskip

In the sequel, whenever it is given an action $A^{\prime }:\left( g,y\right)
\mapsto g\cdot y$ of $G$ on $Y$, we will write $g\cdot x_{0}$, $x_{0}\in
X_{0}$, to indicate the image of the action $A$ on $X_{0}$ defined in (\ref%
{canaction}); a similar convention, based on (\ref{canaction2}), holds for
actions $\underline{A}^{\prime }$ on $\Omega $. Moreover, we will
systematically work under the following assumption.

\bigskip

\textbf{Assumption A -- }Every $Y$-indexed stochastic process $Z$ considered
in the following (not necessarily with a $G$-invariant law) is such that,
for every $x_{0}=\left( \omega ,y\right) \in \Omega \times Y$, the mapping
\begin{equation*}
g\mapsto Z\left( g\cdot x_{0}\right)
\end{equation*}%
is an element of $L^{2}\left( G\right) $.

\bigskip

\textbf{Remark -- }Assumption A can be relaxed in several directions: for
instance, at the cost of some heavier notation, most of the subsequent
results can be immediately extended to stochastic processes $Z$ such that,
for every fixed $y\in Y$, the mapping $g\mapsto Z\left( g\cdot y\right) $ is
in $L^{2}\left( G\right) $ a.s.-$\mathbf{P}$. Note that -- when $G$ is
finite -- Assumption A becomes immaterial.

\bigskip

Now fix an action $A$ of $G$ on $X_{0}$. To every $Y$-indexed stochastic
process $Z$ we associate: the mapping $Z_{\ast }:G\times X_{0}\mapsto
\mathbb{C}$, according to (\ref{jointfunct}), and the mappings $Z_{G}\left[
x_{0}\right] :G\mapsto \mathbb{C}$ and $Z_{X_{0}}\left[ g\right]
:X_{0}\mapsto \mathbb{C}$ as given, respectively, by (\ref{marg1}) for fixed
$x_{0}=\left( \omega ,y\right) \in X_{0}$, and by (\ref{marg2}) for fixed $%
g\in G$. Analogously, for every $\left[ \pi \right] \in \widehat{G}$, the
mapping $Z_{\ast }^{\pi }:G\times X_{0}\mapsto \mathbb{C}$, is defined
according to (\ref{projmarg1}), whereas, for fixed $x_{0}\in X_{0}$ and for
fixed $g\in G$, respectively, $Z_{G}^{\pi }\left[ x_{0}\right] :G\mapsto
\mathbb{C}$ and $Z_{X_{0}}^{\pi }\left[ g\right] :X_{0}\mapsto \mathbb{C}$,
are defined through (\ref{margmarg}). Finally, the mapping $Z^{\pi
}:X_{0}\mapsto \mathbb{C}$ is given by (\ref{projmarg2}).

\begin{proposition}
Under the above notation and assumptions:

\begin{enumerate}
\item for every fixed $x_{0}\in X_{0}$ and for every $\left[ \pi \right] \in
\widehat{G}$, $Z_{G}\left[ x_{0}\right] $ and $Z_{G}^{\pi }\left[ x_{0}%
\right] $ are $\left( G,\mathcal{G}\right) $-measurable functions;

\item for every fixed $g\in G$ and for every $\left[ \pi \right] \in
\widehat{G}$, $Z_{X_{0}}\left[ g\right] $ and $Z_{X_{0}}^{\pi }\left[ g%
\right] $ are $Y$-indexed stochastic processes;

\item if $Z$ has a $G$-invariant law the following three statements hold:
(3-i) for every $g\in G$, the law of $Z_{X_{0}}\left[ g\right] $ is $G$%
-invariant and equal to the law of $Z$; (3-ii) for every $\left[ \pi \right]
\in \widehat{G}$ and $g\in G$, the law of $Z_{X_{0}}^{\pi }\left[ g\right] $
is $G$-invariant and equal to the law of $Z^{\pi }$; (3-iii) the set of
stochastic processes $\left\{ Z,Z^{\pi }:\left[ \pi \right] \in \widehat{G}%
\right\} $ has a $G$-invariant law.
\end{enumerate}
\end{proposition}

\begin{proof}
Points 1.\ and 2. are straightforward. Point (3-i) derives immediately from
the relation: for every $x_{0}\in X_{0}$%
\begin{equation*}
Z_{X_{0}}\left[ g\right] \left( h\cdot x_{0}\right) =Z\left( gh\cdot
x_{0}\right) \text{, \ \ }\forall h\in G,
\end{equation*}%
and the fact that the law of $Z$ is $G$-invariant. To prove point (3-ii), we
can first use the invariance properties of $dg$, as well as the fact that $%
\chi _{\pi }\left( \cdot \right) $ is central for every $\left[ \pi \right]
\in \widehat{G}$, to obtain that for any $h\in G$
\begin{eqnarray}
Z^{\pi }\left( h\cdot x\right) &=&\int_{G}Z\left( gh\cdot x\right) \chi
_{\pi }\left( g^{-1}\right) dg=\int_{G}Z\left( g\cdot x\right) \chi _{\pi
}\left( hg^{-1}\right) dg  \label{proof1} \\
&=&\int_{G}Z\left( g\cdot x\right) \chi _{\pi }\left( g^{-1}h\right)
dg=\int_{G}Z\left( h\cdot \left( g\cdot x\right) \right) \chi _{\pi }\left(
g^{-1}\right) dg  \notag \\
&=&\int_{G}Z_{X}\left[ h\right] \left( g\cdot x\right) \chi _{\pi }\left(
g^{-1}\right) dg,  \notag
\end{eqnarray}%
from which deduce that $Z^{\pi }$ has a $G$-invariant law since, thanks to
point (3-i), $Z_{X}\left[ h\right] $ has the same law as $Z$. To conclude,
just use relation (\ref{projmarg3}), and again point (3-i) applied to the
process $Z^{\pi }$. Point (3-iii) derives immediately from formula (\ref%
{proof1}).
\end{proof}

\bigskip

The following result translates the first part of Theorem 2 into the context
of this section. It shows, in particular, that any $G$-invariant stochastic
process admits a pointwise $L^{2}$-decomposition in terms of simpler $G$%
-invariant stochastic processes, indexed by the elements of $\widehat{G}$.

\begin{theorem}
Let the above notation prevail, and consider an action $A$ of $G$ on $%
X_{0}=\Omega \times Y$. Let $S$ and $Z$ be two $Y$-indexed stochastic
processes verifying Assumption A, and fix $\left[ \pi \right] ,\left[ \sigma %
\right] \in \widehat{G}$ such that $\left[ \pi \right] \neq \left[ \sigma %
\right] $. Then,

\begin{enumerate}
\item for any $\left( \omega _{1},y_{1}\right) ,\left( \omega
_{2},y_{2}\right) \in X_{0}$, $\left\langle S_{G}^{\pi }\left[ \left( \omega
_{1},y_{1}\right) \right] ,Z_{G}^{\sigma }\left[ \left( \omega
_{2},y_{2}\right) \right] \right\rangle _{G}=0$;

\item if, for some $y_{1},y_{2}\in Y$, $S\left( y_{1}\right) ,Z\left(
y_{2}\right) \in L^{2}\left( \mathbf{P}\right) $, then $S^{\pi }\left(
y_{1}\right) ,Z^{\sigma }\left( y_{2}\right) \in L^{2}\left( \mathbf{P}%
\right) $;

\item if the vector $\left( S,Z\right) $ has a $G$-invariant law and $%
S\left( y_{1}\right) ,Z\left( y_{2}\right) \in L^{2}\left( \mathbf{P}\right)
$, then%
\begin{equation}
\mathbf{E}\left[ S^{\pi }\left( y_{1}\right) \overline{Z^{\sigma }\left(
y_{2}\right) }\right] =0;  \label{non corr}
\end{equation}

\item if $S$ has a $G$-invariant law and $S\left( y_{1}\right) \in
L^{2}\left( \mathbf{P}\right) $, then
\begin{equation}
S\left( y_{1}\right) =\sum_{\left[ \pi \right] \in \widehat{G}}S^{\pi
}\left( y_{1}\right) ,  \label{expansion1}
\end{equation}%
where the series on the right hand side is orthogonal and convergent in $%
L^{2}\left( \mathbf{P}\right) $.
\end{enumerate}
\end{theorem}

\begin{proof}
Point 1. is a direct consequence of Theorem 2-1, whereas point 2. derives
from the inequality%
\begin{eqnarray*}
\mathbf{E}\left[ \left\vert S^{\pi }\left( y_{1}\right) \right\vert ^{2}%
\right] &=&d_{\pi }^{2}\int_{\Omega }\mathbf{P}\left( d\omega \right)
\left\vert \int_{G}\chi _{\pi }\left( g\right) S\left( g^{-1}\cdot \left(
\omega ,y_{1}\right) \right) dg\right\vert ^{2} \\
&\leq &d_{\pi }^{2}\alpha _{\pi }^{2}\int_{\Omega }\mathbf{P}\left( d\omega
\right) \int_{G}\left\vert S\left( g^{-1}\cdot \left( \omega ,y_{1}\right)
\right) \right\vert ^{2}dg \\
&=&d_{\pi }^{2}\alpha _{\pi }^{2}\int_{G}\mathbf{E}\left[ \left\vert S\left(
g^{-1}\cdot y_{1}\right) \right\vert ^{2}\right] dg=d_{\pi }^{2}\alpha _{\pi
}^{2}\mathbf{E}\left[ \left\vert S\left( y_{1}\right) \right\vert ^{2}\right]
<+\infty \text{,}
\end{eqnarray*}%
and a similar calculation for $Z^{\pi }$. To see point 3., just write, due
to the $G$-invariance of $\left( S,Z\right) $ and the fact that $%
\int_{G}dg=1 $,%
\begin{eqnarray*}
\mathbf{E}\left[ S^{\pi }\left( y_{1}\right) \overline{Z^{\sigma }\left(
y_{2}\right) }\right] &=&\int_{G}\mathbf{E}\left[ S^{\pi }\left( g\cdot
y_{1}\right) \overline{Z^{\sigma }\left( g\cdot y_{2}\right) }\right] dg \\
&=&\int_{\Omega }\left\{ \int_{G}S^{\pi }\left[ \left( \omega ,y_{1}\right) %
\right] \left( g\right) \overline{Z^{\sigma }\left[ \left( \omega
,y_{2}\right) \right] \left( g\right) }dg\right\} \mathbf{P}\left( d\omega
\right) \\
&=&\int_{\Omega }\left\langle S_{G}^{\pi }\left[ \left( \omega ,y_{1}\right) %
\right] ,Z_{G}^{\sigma }\left[ \left( \omega ,y_{2}\right) \right]
\right\rangle _{G}\mathbf{P}\left( d\omega \right) =0,
\end{eqnarray*}%
where we have used a Fubini theorem, as well as point 1., with $\omega
_{1}=\omega _{2}=\omega .$ To prove point 4., let $\left[ \pi \left(
i\right) \right] $, $i=1,2,...$, be an enumeration of the elements of $%
\widehat{G}$, and observe that, according to (\ref{decomposition1}), for
every $x_{0}=\left( \omega ,y_{1}\right) \in X_{0}$

\begin{equation*}
\lim_{N\rightarrow +\infty }\int_{G}\left\vert S_{G}\left[ x_{0}\right]
\left( g\right) -\sum_{i=1}^{N}S_{G}^{\pi \left( i\right) }\left[ x_{0}%
\right] \left( g\right) \right\vert ^{2}dg=0,
\end{equation*}%
and also, thanks to (\ref{fouriercoeff1}),
\begin{equation*}
\int_{G}\left\vert S_{G}\left[ x_{0}\right] \left( g\right)
-\sum_{i=1}^{N}S_{G}^{\pi \left( i\right) }\left[ x_{0}\right] \left(
g\right) \right\vert ^{2}dg=\int_{G}\left\vert \sum_{i=N+1}^{\infty
}S_{G}^{\pi \left( i\right) }\left[ x_{0}\right] \left( g\right) \right\vert
^{2}dg\leq \int_{G}\left\vert S_{G}\left[ x_{0}\right] \left( g\right)
\right\vert ^{2}dg.
\end{equation*}%
Now observe that, for fixed $y_{1}\in Y$, the random variable $\omega
\mapsto \int_{G}\left\vert S\left( g\cdot \left( \omega ,y_{1}\right)
\right) \right\vert ^{2}dg$ is in $L^{1}\left( \mathbf{P}\right) $, since,
due to the $G$-invariance of the law of $S$,%
\begin{equation*}
\int_{\Omega }\left[ \int_{G}\left\vert S\left( g\cdot \left( \omega
,y_{1}\right) \right) \right\vert ^{2}dg\right] \mathbf{P}\left( d\omega
\right) =\int_{G}\mathbf{E}\left[ \left\vert S\left( g\cdot y_{1}\right)
\right\vert ^{2}\right] dg=\mathbf{E}\left[ \left\vert S\left( y_{1}\right)
\right\vert ^{2}\right] <+\infty .
\end{equation*}%
Finally, according to Proposition 3-3-iii, the class $\left\{ S,S^{\pi }:%
\left[ \pi \right] \in \widehat{G}\right\} $ has a $G$-invariant law, and
therefore%
\begin{eqnarray*}
\mathbf{E}\left[ \left\vert S\left( y_{1}\right) -\sum_{i=1}^{N}S^{\pi
\left( i\right) }\left( y_{1}\right) \right\vert ^{2}\right] &=&\int_{G}%
\mathbf{E}\left[ \left\vert S\left( g\cdot y_{1}\right)
-\sum_{i=1}^{N}S^{\pi \left( i\right) }\left( g\cdot y_{1}\right)
\right\vert ^{2}\right] dg \\
&=&\mathbf{E}\left[ \int_{G}\left\vert S_{G}\left[ x_{0}\right] \left(
g\right) -\sum_{i=1}^{N}S_{G}^{\pi \left( i\right) }\left[ x_{0}\right]
\left( g\right) \right\vert ^{2}dg\right] \underset{N\rightarrow +\infty }{%
\rightarrow }0,
\end{eqnarray*}%
due to an application of the dominated convergence theorem.
\end{proof}

\bigskip

When $G$ is finite, formula (\ref{expansion1}) holds even if the law of $S$
is not $G$-invariant (but, in this case, the series is not necessarily
orthogonal in $L^{2}\left( \mathbf{P}\right) $). We now apply Theorem 2 to
further characterize actions of the specific form (\ref{canaction}). Observe
that the following Theorem applies to processes whose laws are not
necessarily $G$-invariant.

\bigskip

\begin{theorem}
Let the action $A:G\times X_{0}\mapsto X_{0}$ be such that, $\forall \left(
\omega ,y\right) \in X_{0}$, $A\left( g,\left( \omega ,y\right) \right)
=\left( \omega ,A^{\prime }\left( g,y\right) \right) $, where $A^{\prime }$
is an action on $\left( Y,\mathcal{Y}\right) $. Consider moreover two $Y$%
-indexed stochastic processes $S,Z$ (not necessarily with $G$-invariant
laws), as well as a $\sigma $-finite positive measure $\mu $ on $\left( Y,%
\mathcal{Y}\right) $, which is invariant with respect to the action $%
A^{\prime }$ of $G$ on $Y$. Suppose that, for every fixed $\omega ^{\ast
}\in \Omega $, the applications $\left( g,y\right) \mapsto Z\left( \omega
^{\ast },A^{\prime }\left( g,y\right) \right) $ and $\left( g,y\right)
\mapsto S\left( \omega ^{\ast },A^{\prime }\left( g,y\right) \right) $ are
elements of $L^{2}\left( dg\times \mu \left( dy\right) \right) $, and also
that, for every fixed $\left( \omega ^{\ast },g^{\ast }\right) \in \Omega
\times G$, the mappings $y\mapsto Z\left( \omega ^{\ast },A^{\prime }\left(
g^{\ast },y\right) \right) $ and $y\mapsto S\left( \omega ^{\ast },A^{\prime
}\left( g^{\ast },y\right) \right) $ are in $L^{2}\left( \mu \left(
dy\right) \right) $. Then,

\begin{enumerate}
\item for every fixed $\left( \omega ^{\ast },g^{\ast }\right) \in \Omega
\times G$, and for every $\left[ \pi \right] ,\left[ \sigma \right] \in
\widehat{G}$ such that $\left[ \pi \right] \neq \left[ \sigma \right] $,
\begin{equation}
\int_{Y}S_{X_{0}}^{\pi }\left[ g^{\ast }\right] \left( \omega ^{\ast
},y\right) \overline{Z_{X_{0}}^{\sigma }\left[ g^{\ast }\right] \left(
\omega ^{\ast },y\right) }\mu \left( dy\right) =0;  \label{traj-ortho}
\end{equation}

\item for every fixed $\omega ^{\ast }\in \Omega $,
\begin{equation}
S\left( \omega ^{\ast },y\right) =\sum_{\left[ \pi \right] \in \widehat{G}%
}S^{\pi }\left( \omega ^{\ast },y\right) \text{ \ \ and \ \ }Z\left( \omega
^{\ast },y\right) =\sum_{\left[ \pi \right] \in \widehat{G}}Z^{\pi }\left(
\omega ^{\ast },y\right) ,  \label{traj-decomposition}
\end{equation}%
where the two series are orthogonal and convergent in $L^{2}\left( \mu
\left( dy\right) \right) $, and therefore%
\begin{equation}
\left\langle S\left( \omega ^{\ast },\cdot \right) ,Z\left( \omega ^{\ast
},\cdot \right) \right\rangle _{L^{2}\left( \mu \left( dy\right) \right)
}=\sum_{\left[ \pi \right] \in \widehat{G}}\left\langle S^{\pi }\left(
\omega ^{\ast },\cdot \right) ,Z^{\pi }\left( \omega ^{\ast },\cdot \right)
\right\rangle _{L^{2}\left( \mu \left( dy\right) \right) };  \label{traj-PP}
\end{equation}

\item if moreover $Z\left( \omega ,y\right) \in L^{2}\left( \mathbf{P}\left(
d\omega \right) \times \mu \left( dy\right) \right) $, then,
\begin{equation*}
Z^{\pi }\left( \omega ,y\right) \in L^{2}\left( \mathbf{P}\left( d\omega
\right) \times \mu \left( dy\right) \right) \text{, for every }\left[ \pi %
\right] \in \widehat{G}\text{,}
\end{equation*}%
and%
\begin{equation}
Z\left( \omega ,y\right) =\sum_{\left[ \pi \right] \in \widehat{G}}Z^{\pi
}\left( \omega ,y\right) ,  \label{traj-series}
\end{equation}%
where the series is orthogonal and convergent in $L^{2}\left( \mathbf{P}%
\left( d\omega \right) \times \mu \left( dy\right) \right) $.
\end{enumerate}
\end{theorem}

\begin{proof}
(\textbf{1.}) For every $\omega ^{\ast }\in \Omega $, the measure $\nu
^{\ast }$ on $\left( X_{0},\mathcal{X}_{0}\right) =\left( \Omega \times Y,%
\mathcal{F}\otimes \mathcal{Y}\right) $, defined by $\nu ^{\ast }\left(
d\omega ,dy\right) =\delta _{\omega ^{\ast }}\left( d\omega \right) \times
\mu \left( dy\right) $, where $\delta _{\omega ^{\ast }}$ is the Dirac mass
at $\omega ^{\ast }$, is invariant with respect to the action $A$ of $G$ on $%
X_{0}$. Moreover, it is easily seen that the assumptions in the statement
imply that $\nu ^{\ast }$ satisfies all the hypotheses of Theorem 2-3, so
that formula (\ref{traj-ortho}) follows immediately, by observing that, for
every $g\in G$,%
\begin{equation*}
\int_{Y}S_{X_{0}}^{\pi }\left[ g\right] \left( \omega ^{\ast },y\right)
\overline{Z_{X_{0}}^{\sigma }\left[ g\right] \left( \omega ^{\ast },y\right)
}\mu \left( dy\right) =\int_{X_{0}}S_{X_{0}}^{\pi }\left[ g\right] \left(
x_{0}\right) \overline{Z_{X_{0}}^{\sigma }\left[ g\right] \left(
x_{0}\right) }\nu ^{\ast }\left( dx_{0}\right) .
\end{equation*}

(\textbf{2.}) This is a direct consequence of Theorem 2-5 (in the case $g=e$%
).

(\textbf{3}.) First observe that $Z^{\pi }\left( \omega ,y\right) \in
L^{2}\left( \mathbf{P}\left( d\omega \right) \times \mu \left( dy\right)
\right) $, since, thanks to the $G$-invariance of $\mu $,%
\begin{eqnarray*}
\mathbf{E}\left[ \int_{Y}\left\vert Z^{\pi }\left( y\right) \right\vert
^{2}\mu \left( dy\right) \right] &\leq &d_{\pi }^{2}\alpha _{\pi
}^{2}\int_{Y}\mu \left( dy\right) \int_{\Omega }\mathbf{P}\left( d\omega
\right) \int_{G}dg\left\vert Z\left( \omega ,g^{-1}\cdot y\right)
\right\vert ^{2} \\
&=&d_{\pi }^{2}\alpha _{\pi }^{2}\int_{G}\int_{Y}\mathbf{E}\left[ \left\vert
Z\left( g^{-1}\cdot y\right) \right\vert ^{2}\right] \mu \left( dy\right) dg
\\
&=&d_{\pi }^{2}\alpha _{\pi }^{2}\int_{Y}\mathbf{E}\left[ \left\vert Z\left(
y\right) \right\vert ^{2}\right] \mu \left( dy\right) <+\infty \text{.}
\end{eqnarray*}%
The rest of the proof is similar to the proof of Theorem 3-4. Let indeed $%
\left[ \pi \left( i\right) \right] $, $i=1,2,...$, be an enumeration of $%
\widehat{G}$, and observe that Theorem 2-4 (formula (\ref{twovardecomp}))
implies that, for every $\omega ^{\ast }\in \Omega $,%
\begin{eqnarray*}
&&\lim_{N\rightarrow +\infty }\int_{Y}\int_{G}\left\vert Z\left( \omega
^{\ast },g\cdot y\right) -\sum_{i=1}^{N}Z^{\pi \left( i\right) }\left(
\omega ^{\ast },g\cdot y\right) \right\vert ^{2}dg\mu \left( dy\right) \\
&=&\lim_{N\rightarrow +\infty }\int_{X_{0}}\int_{G}\left\vert Z\left( g\cdot
x_{0}\right) -\sum_{i=1}^{N}Z^{\pi \left( i\right) }\left( g\cdot
x_{0}\right) \right\vert ^{2}dg\nu ^{\ast }\left( dx_{0}\right) =0,
\end{eqnarray*}%
and (\ref{fourcoeff2}) yields also%
\begin{eqnarray*}
\int_{Y}\int_{G}\left\vert Z\left( \omega ^{\ast },g\cdot y\right)
-\sum_{i=1}^{N}Z^{\pi \left( i\right) }\left( \omega ^{\ast },g\cdot
y\right) \right\vert ^{2}dg\mu \left( dy\right) &\leq
&\int_{Y}\int_{G}\left\vert Z\left( g\cdot y\right) \right\vert ^{2}dg\mu
\left( dy\right) \\
&=&\int_{Y}\left\vert Z\left( y\right) \right\vert ^{2}\mu \left( dy\right)
\in L^{1}\left( \mathbf{P}\right) \text{,}
\end{eqnarray*}%
since $Z\in L^{2}\left( \mathbf{P}\left( d\omega \right) \times \mu \left(
dy\right) \right) $ by assumption. Since $\mu $ is $G$-invariant, and by
dominated convergence,%
\begin{eqnarray*}
\mathbf{E}\left[ \int_{Y}\left\vert Z\left( y\right) -\sum_{i=1}^{N}Z^{\pi
\left( i\right) }\left( y\right) \right\vert ^{2}\mu \left( dy\right) \right]
&=&\int_{G}\mathbf{E}\left[ \int_{Y}\left\vert Z\left( g\cdot y\right)
-\sum_{i=1}^{N}Z^{\pi \left( i\right) }\left( g\cdot y\right) \right\vert
^{2}\mu \left( dy\right) \right] dg \\
&=&\mathbf{E}\left[ \int_{Y}\int_{G}\left\vert Z\left( g\cdot y\right)
-\sum_{i=1}^{N}Z^{\pi \left( i\right) }\left( g\cdot y\right) \right\vert
^{2}dg\mu \left( dy\right) \right] \underset{N\rightarrow +\infty }{%
\rightarrow }0.
\end{eqnarray*}
\end{proof}

\subsection{Gaussian processes}

Keep the previous notation and assumptions (in particular, Assumption A
holds throughout the following). In this paragraph, we apply the above
established results to the case of a \textbf{two-dimensional} \textbf{%
real-valued Gaussian process }of the type\textbf{\ }%
\begin{equation*}
\left( Z_{1},Z_{2}\right) :\Omega \times Y\mapsto \mathbb{R}^{2}:\left(
\omega ,y\right) \mapsto \left( Z_{1}\left( \omega ,y\right) ,Z_{2}\left(
\omega ,y\right) \right)
\end{equation*}%
with a \textbf{covariance structure} given by
\begin{equation}
R_{i,j}\left( y_{1},y_{2}\right) =\mathbf{E}\left[ Z_{i}\left( y_{1}\right)
Z_{j}\left( y_{2}\right) \right] \text{, \ \ }i,j=1,2\text{, \ \ }%
y_{1},y_{2}\in Y\text{.}  \label{COV}
\end{equation}

Note that our definition of two-dimensional Gaussian process also covers the
case $Z_{1}=Z_{2}$. In this paragraph, we will consider exclusively actions
of the type (\ref{canaction}), where $A^{\prime }$ is an action of the
topological compact group $G$ on $Y$. Note that, under such assumptions, $%
\left( Z_{1},Z_{2}\right) $ has a $G$-invariant law if, and only if,
\begin{equation}
R_{i,j}\left( g\cdot y_{1},g\cdot y_{2}\right) =R_{i,j}\left(
y_{1},y_{2}\right) \text{, \ for every }g\in G\text{,\ }i,j=1,2\text{, and }%
y_{1},y_{2}\in Y.  \label{ICOV}
\end{equation}

When the function $R_{i,j}$ satisfies (\ref{ICOV}), we say that $R_{i,j}$ is
a $G$-\textbf{invariant covariance function}.

\bigskip

In the sequel, the Cartesian product\textbf{\ }$G\times G=G^{2}$ is
systematically endowed with the \textbf{product group }structure, as
described e.g. in \cite[Section 3.2]{Serre}. The generic element of $G^{2}$
is noted $\left( g_{1},g_{2}\right) $; $G^{2}$ is again a topological and
compact group, with Haar measure given by $dg_{1}\times dg_{2}$. Recall (see
again \cite[Section 1.5 and 3.2]{Serre}) that $\left[ \rho \right] \in
\widehat{G^{2}}$ if, and only if, $\left[ \rho \right] =\left[ \pi _{1}%
\right] \otimes \left[ \pi _{2}\right] $, where $\left( \left[ \pi _{1}%
\right] ,\left[ \pi _{2}\right] \right) \in \widehat{G}\times \widehat{G}$,
and $\otimes $ stands for the (tensor) product between representations. The
following assumption will hold for the rest of the section.

\bigskip

\textbf{Assumption B -- }For every two-dimensional Gaussian process $\left(
Z_{1},Z_{2}\right) $ considered in the sequel, and for every fixed $%
y_{1},y_{2}\in Y$, the application
\begin{equation}
\left( R_{i,j}\right) _{G^{2}}\left[ y_{1},y_{2}\right] :G\times G\mapsto
\mathbb{R}:\left( g_{1},g_{2}\right) \mapsto R_{i,j}\left( g_{1}\cdot
y_{1},g_{2}\cdot y_{2}\right)  \label{margcov}
\end{equation}%
(see (\ref{ICOV}), and observe that (\ref{margcov}) is consistent with the
notation introduced in (\ref{marg1})) is an element of $L^{2}\left(
G^{2}\right) $, for every $i,j=1,2.$

\bigskip

Again, if $G$ is finite, Assumption B is redundant. Given $\widehat{G^{2}}%
\ni \left[ \rho \right] =\left[ \pi _{1}\right] \otimes \left[ \pi _{2}%
\right] $, we define, for fixed $y_{1},y_{2}\in Y$,
\begin{equation}
\left( R_{i,j}\right) _{G^{2}}^{\rho }\left[ y_{1},y_{2}\right] =\left(
R_{i,j}\right) _{G^{2}}^{\pi _{1}\otimes \pi _{2}}\left[ y_{1},y_{2}\right]
\label{cov-projection}
\end{equation}%
according to (\ref{margmarg}). The following result, which is a consequence
of Theorem 4 will lead to a very general version of Watson's duplication
identity.

\bigskip

\begin{proposition}
Let $\left( Z_{1},Z_{2}\right) $ be a two dimensional real-valued Gaussian
process with a $G$-invariant law. Then,

\begin{enumerate}
\item the collection of (possibly complex-valued) stochastic processes $%
\left\{ Z_{1},Z_{2},Z_{1}^{\pi },Z_{2}^{\sigma }:\left[ \pi \right] ,\left[
\sigma \right] \in \widehat{G}\right\} $ is jointly Gaussian;

\item for every $\left[ \pi \right] ,\left[ \sigma \right] \in \widehat{G}$
such that $\left[ \pi \right] \neq \left[ \sigma \right] $ and $\chi _{\pi }$
is real valued, the two processes $Z_{i}^{\pi }$ and $Z_{j}^{\sigma }$ are
independent for every $i,j=1,2$;

\item for every $\left[ \pi \right] ,\left[ \sigma \right] $ as at point 2.,
$\left( R_{i,j}\right) _{G^{2}}^{\pi \otimes \sigma }\left[ y_{1},y_{2}%
\right] =0$, for every $i,j=1,2;$

\item for $i=1,2$, $Z_{i}\left( y\right) =\sum_{\left[ \pi \right] \in
\widehat{G}}Z_{i}^{\pi }\left( y\right) $ in $L^{2}\left( \mathbf{P}\right) $%
.
\end{enumerate}
\end{proposition}

\begin{proof}
Point 1. is immediate, since the action $A$ of $G$ on $X_{0}$ has the form (%
\ref{canaction}). Since $\chi _{\pi }$ is real-valued, $Z_{i}^{\pi }$ is
also real-valued, and moreover, for every $y_{1},y_{2}\in Y$, according to
Theorem 4-3, $\mathbf{E}\left[ Z_{i}^{\pi }\left( y_{1}\right) \overline{%
Z_{j}^{\sigma }\left( y_{2}\right) }\right] =0$, thus implying that $%
Z_{i}^{\pi }\left( y_{1}\right) $ is independent of both the real and
imaginary parts of $Z_{j}^{\sigma }\left( y_{2}\right) $. This concludes the
proof of point 2. To see point 3., just write, for $h_{1},h_{2}\in G$%
\begin{eqnarray*}
\left( R_{i,j}\right) _{G^{2}}^{\pi \otimes \sigma }\left[ y_{1},y_{2}\right]
\left( h_{1},h_{2}\right) &=&\int_{G}dg_{1}\int_{G}dg_{2}\chi _{\pi }\left(
g_{1}\right) \chi _{\sigma }\left( g_{2}\right) \left( R_{i,j}\right)
_{G^{2}}\left[ y_{1},y_{2}\right] \left(
g_{1}^{-1}h_{1},g_{2}^{-1}h_{2}\right) \\
&=&\mathbf{E}\left[ Z_{i}^{\pi }\left( h_{1}\cdot y_{1}\right) Z_{j}^{\sigma
}\left( h_{2}\cdot y_{2}\right) \right] =0\text{,}
\end{eqnarray*}%
due to point 2.. Point 4. comes immediately from Theorem 4-4.
\end{proof}

\bigskip

Of course, point 1. of Proposition 3 still holds when the law of the
Gaussian process $\left( Z_{1},Z_{2}\right) $ is not $G$-invariant. The
combination of Theorem 5 and Proposition 6 yields immediately the following

\bigskip

\begin{proposition}
Let $G$ be such that $\chi _{\pi }$ is real-valued for every $\left[ \pi %
\right] \in \widehat{G}$. Let $\left( Z_{1},Z_{2}\right) $ be a two
dimensional real-valued Gaussian process with a $G$-invariant law, and
consider a $G$-invariant, $\sigma $-finite and positive measure $\mu $ on $%
\left( Y,\mathcal{Y}\right) $. Suppose that, for any fixed $\omega ^{\ast
}\in \Omega $ and $i=1,2$, the mapping $\left( g,y\right) \mapsto
Z_{i}\left( \omega ^{\ast },g\cdot y\right) $ is in $L^{2}\left( dg\times
\mu \left( dy\right) \right) $, and also that, for every fixed $\left(
\omega ^{\ast },g^{\ast }\right) \in \Omega \times G$, the mapping $y\mapsto
Z_{i}\left( \omega ^{\ast },g^{\ast }\cdot y\right) $ is an element of $%
L^{2}\left( \mu \left( dy\right) \right) $. Then, for every $i,j=1,2$,

\begin{enumerate}
\item the Gaussian processes $Z_{i}^{\pi }\left( \omega ,y\right) $ and $%
Z_{j}^{\sigma }\left( \omega ,y\right) $ are independent for every $\left[
\pi \right] \neq \left[ \sigma \right] $, and orthogonal in $L^{2}\left( \mu
\left( dy\right) \right) $ for every $\omega \in \Omega $;

\item for every $\left[ \pi \right] \in \widehat{G}$,%
\begin{equation}
\mathbf{E}\left[ Z_{i}^{\pi }\left( y_{1}\right) Z_{j}^{\pi }\left(
y_{2}\right) \right] =\left( R_{i,j}\right) _{G^{2}}^{\pi \otimes \pi }\left[
y_{1},y_{2}\right] \left( e,e\right) =R_{i,j}^{\pi \otimes \pi }\left(
y_{1},y_{2}\right) \text{;}  \label{reprCov}
\end{equation}

\item $Z_{i}\left( \omega ,y\right) =\sum_{\left[ \pi \right] \in \widehat{G}%
}Z_{i}^{\pi }\left( \omega ,y\right) $ both in $L^{2}\left( \mu \left(
dy\right) \right) $ for every $\omega \in \Omega $ and in $L^{2}\left(
\mathbf{P}\left( d\omega \right) \times \mu \left( dy\right) \right) $;

\item for every $\lambda \in \mathbb{R}$,
\begin{equation}
\mathbf{E}\left[ \exp \left( \mathtt{i}\lambda \int_{Y}Z_{i}\left( y\right)
Z_{j}\left( y\right) \mu \left( dy\right) \right) \right] =\prod_{\left[ \pi %
\right] \in \widehat{G}}\mathbf{E}\left[ \exp \left( \mathtt{i}\lambda
\int_{Y}Z_{i}^{\pi }\left( y\right) Z_{j}^{\pi }\left( y\right) \mu \left(
dy\right) \right) \right]  \label{GenWats}
\end{equation}
\end{enumerate}
\end{proposition}

\bigskip

\textbf{Example }(A group-theoretic proof of the (polarized) Watson's
identity) -- As a first illustration of our techniques, we shall obtain a
class of identities in law -- between functionals of two correlated Brownian
bridges -- extending Watson's identity (\ref{Watson1961}). Our method of
proof, which is directly based on the discussion contained in this
paragraph, generalizes the simple proof of (\ref{Watson1961}) given by the
second author in \cite{Pycke}, and will motivate the content of the
subsequent section. To this end, we consider a two-dimensional Brownian
bridge $\underline{b}=\left\{ b_{1}\left( t\right) ,b_{2}\left( t\right)
:t\in \left[ 0,1\right] \right\} $ with correlation parameter equal to $\rho
\in \left[ 0,1\right] $. This means that $\underline{b}$ is a
two-dimensional, real-valued Gaussian process such that, for every $s,t\in %
\left[ 0,1\right] $, $\mathbf{E}\left[ b_{i}\left( s\right) b_{i}\left(
t\right) \right] =s\wedge t-st$, $i=1,2$, and $\mathbf{E}\left[ b_{1}\left(
s\right) b_{2}\left( t\right) \right] =\rho \times \left( s\wedge
t-st\right) $. By $\underline{b}_{\ast }=\left\{ b_{\ast 1}\left( t\right)
,b_{\ast 2}\left( t\right) :t\in \left[ 0,1\right] \right\} $, we denote an
independent copy of $\underline{b}$, and we also write, for $i=1,2$ and $%
t\in \left[ 0,1\right] $,%
\begin{equation}
v_{i}\left( t\right) =b_{i}\left( t\right) -\int_{0}^{1}b_{i}\left( s\right)
ds\text{ \ \ and \ \ }v_{\ast i}\left( t\right) =b_{\ast i}\left( t\right)
-\int_{0}^{1}b_{\ast i}\left( s\right) ds.  \label{compensatedBr}
\end{equation}%
Now consider the group $G=\left\{ e,g\right\} \simeq \mathbb{Z}/2\mathbb{Z}$%
, where $e$ stands again for the identity element. It is plain (see e.g.
\cite[Chapter 2]{Serre}) that in this case $\widehat{G}=\left\{ \left[ \pi
_{u}\right] ,\left[ \pi _{a}\right] \right\} $, where $\left[ \pi _{u}\right]
$ and $\left[ \pi _{a}\right] $ are the equivalence classes, respectively of
the \textit{unity} and of the \textit{alternating} representation; in
particular, $\chi _{\pi _{u}}\left( e\right) =\chi _{\pi _{u}}\left(
g\right) =1$, and $\chi _{\pi _{a}}\left( e\right) =1=-\chi _{\pi
_{a}}\left( g\right) $. We fix the following elementary action of $G$ on $%
\left[ 0,1\right] $: $e\cdot t=t$ and $g\cdot t=1-t$, $\forall t\in \left[
0,1\right] $. It is well known that $\underline{b}$, and therefore the
vector $\left( v_{1},v_{2}\right) $, has a $G$-invariant law, so that the
content of Proposition 6 can be directly applied. To do this, we first set,
according to (\ref{projmarg2}) and for $i=1,2$ and $t\in \left[ 0,1\right] $%
, and since $d_{\pi _{u}}=d_{\pi _{a}}=1$,%
\begin{eqnarray*}
v_{i}^{\pi _{u}}\left( t\right) &=&\frac{d_{\pi _{u}}}{\left\vert
G\right\vert }\left\{ \chi _{\pi _{u}}\left( e\right) v_{i}\left(
e^{-1}\cdot t\right) +\chi _{\pi _{u}}\left( g\right) v_{i}\left(
g^{-1}\cdot x\right) \right\} =\frac{1}{2}\left( b_{i}\left( t\right)
+b_{i}\left( 1-t\right) \right) -\int_{0}^{1}b_{i}\left( s\right) ds \\
&=&\frac{1}{2}\left( b_{i}\left( t\right) +b_{i}\left( 1-t\right) \right) -%
\frac{1}{2}\int_{0}^{1}\left( b_{i}\left( s\right) +b_{i}\left( 1-s\right)
\right) ds \\
v_{i}^{\pi _{a}}\left( t\right) &=&\frac{d_{\pi _{a}}}{\left\vert
G\right\vert }\left\{ \chi _{\pi _{a}}\left( e\right) v_{i}\left(
e^{-1}\cdot t\right) +\chi _{\pi a}\left( g\right) v_{i}\left( g^{-1}\cdot
x\right) \right\} =\frac{1}{2}\left( b_{i}\left( t\right) -b_{i}\left(
1-t\right) \right) \text{,}
\end{eqnarray*}%
and an analogous definition holds for $v_{\ast i}^{\pi _{u}}$ and $v_{\ast
i}^{\pi _{a}}$, $i=1,2$. Now observe that Proposition 6-2 (in the case $%
\left( Z_{1},Z_{2}\right) =\left( v_{1},v_{2}\right) $) implies that, for
any $i,j=1,2$, the two processes $v_{i}^{\pi _{u}}$ and $v_{j}^{\pi _{a}}$
are independent. Moreover, the restriction of Lebesgue measure to $\left[ 0,1%
\right] $ is trivially $G$-invariant, so that all assumptions of Proposition
7 are satisfied (again with $\left( Z_{1},Z_{2}\right) =\left(
v_{1},v_{2}\right) $ and $\mu $ equal to Lebesgue measure) and therefore%
\begin{eqnarray*}
&&\int_{0}^{1}v_{1}\left( t\right) v_{2}\left( t\right) dt\overset{law}{=}%
\int_{0}^{1}v_{1}^{\pi _{u}}\left( t\right) v_{2}^{\pi _{u}}\left( t\right)
dt+\int_{0}^{1}v_{\ast 1}^{\pi _{a}}\left( t\right) v_{\ast 2}^{\pi
_{a}}\left( t\right) dt \\
&=&\frac{1}{4}\int_{0}^{1}\left( b_{1}\left( t\right) +b_{1}\left(
1-t\right) -\int_{0}^{1}\left( b_{1}\left( s\right) +b_{1}\left( 1-s\right)
\right) ds\right) \times \\
&&\times \left( b_{2}\left( t\right) +b_{2}\left( 1-t\right)
-\int_{0}^{1}\left( b_{2}\left( s\right) +b_{2}\left( 1-s\right) \right)
ds\right) dt \\
&&+\frac{1}{4}\int_{0}^{1}\left( b_{\ast 1}\left( t\right) -b_{\ast 1}\left(
1-t\right) \right) \left( b_{\ast 2}\left( t\right) -b_{\ast 2}\left(
1-t\right) \right) dt.
\end{eqnarray*}%
Next, consider a correlated two-dimensional standard Brownian motion $%
\underline{W}=\left\{ W_{1}\left( t\right) ,W_{2}\left( t\right) :t\in \left[
0,1\right] \right\} $ with correlation $\rho $,\footnote{%
That is, $\underline{W}$ is a two-dimensional Gaussian process such that,
for $i=1,2$ and $s,t\in \left[ 0,1\right] $, $\mathbf{E}\left[ W_{i}\left(
s\right) W_{i}\left( t\right) \right] =s\wedge t$ and $\mathbf{E}\left[
W_{1}\left( s\right) W_{2}\left( t\right) \right] =\rho \times \left(
s\wedge t\right) $.} and independent of $\underline{b}$. Routine
computations show the following identities in law:
\begin{eqnarray*}
&&\left\{ b_{1}\left( t\right) +b_{1}\left( 1-t\right) ,b_{2}\left( t\right)
+b_{2}\left( 1-t\right) :t\in \left[ 0,1/2\right] \right\} \overset{law}{=}%
\left\{ W_{1}\left( 2t\right) ,W_{2}\left( 2t\right) :t\in \left[ 0,1/2%
\right] \right\} \\
&&\left\{ b_{1}\left( t\right) -b_{1}\left( 1-t\right) ,b_{2}\left( t\right)
-b_{2}\left( 1-t\right) :t\in \left[ 0,1/2\right] \right\} \overset{law}{=}%
\left\{ b_{1}\left( 2t\right) ,b_{2}\left( 2t\right) :t\in \left[ 0,1/2%
\right] \right\} ,
\end{eqnarray*}%
implying that%
\begin{equation*}
\int_{0}^{1}v_{1}\left( t\right) v_{2}\left( t\right) \overset{law}{=}\frac{1%
}{4}\int_{0}^{1}\left( W_{1}\left( t\right) -\int_{0}^{1}W_{1}\left(
s\right) ds\right) \left( W_{2}\left( t\right) -\int_{0}^{1}W_{2}\left(
s\right) ds\right) dt+\frac{1}{4}\int_{0}^{1}b_{1}\left( t\right)
b_{2}\left( t\right) dt.
\end{equation*}%
We eventually use some standard arguments (see e.g. \cite[Lemma 2]{Chou}) to
prove that
\begin{equation*}
\int_{0}^{1}\left( W_{1}\left( t\right) -\int_{0}^{1}W_{1}\left( s\right)
ds\right) \left( W_{2}\left( t\right) -\int_{0}^{1}W_{2}\left( s\right)
ds\right) dt\overset{law}{=}\int_{0}^{1}b_{1}\left( t\right) b_{2}\left(
t\right) dt
\end{equation*}%
and therefore
\begin{equation*}
\int_{0}^{1}\left( b_{1}\left( t\right) -\int_{0}^{1}b_{1}\left( s\right)
ds\right) \left( b_{2}\left( t\right) -\int_{0}^{1}b_{2}\left( s\right)
ds\right) dt\overset{law}{=}\frac{1}{4}\int_{0}^{1}\left( b_{1}\left(
t\right) b_{2}\left( t\right) +b_{\ast 1}\left( t\right) b_{\ast 2}\left(
t\right) \right) dt
\end{equation*}%
(Watson's identity (\ref{Watson1961}) can be obtained by setting $\rho =1$).

\bigskip

\textbf{Remark -- }By using e.g. \cite[Proposition 2]{Chou}, we obtain that,
for $\lambda >0$ sufficiently small and $\rho \in \left[ 0,1\right] $%
\begin{eqnarray*}
&&\mathbf{E}\left[ \exp \left( \lambda ^{2}\int_{0}^{1}\left( b_{1}\left(
t\right) -\int_{0}^{1}b_{1}\left( s\right) ds\right) \left( b_{2}\left(
t\right) -\int_{0}^{1}b_{2}\left( s\right) ds\right) dt\right) \right] \\
&=&\frac{\left( \lambda /2\right) ^{2}\sqrt{1-\rho ^{2}}}{\sin \frac{\lambda
}{2}\sqrt{1+\rho }\sinh \frac{\lambda }{2}\sqrt{1-\rho }}
\end{eqnarray*}

\bigskip

Note that the $G$-invariant process $\left( v_{1},v_{2}\right) $, introduced
in formula (\ref{compensatedBr}) of the previous example, has the remarkable
property that
\begin{equation}
\int_{0}^{1}v_{1}^{\pi _{u}}\left( t\right) v_{2}^{\pi _{u}}\left( t\right)
dt\overset{law}{=}\int_{0}^{1}v_{\ast 1}^{\pi _{a}}\left( t\right) v_{\ast
2}^{\pi _{a}}\left( t\right) dt.  \label{WatRel}
\end{equation}

In the next paragraph we shall establish necessary and sufficient conditions
to ensure that, in the case of a finite $G$, a $G$-invariant Gaussian
process $\left( Z_{1},Z_{2}\right) $ (with some special structure) is such
that
\begin{equation}
\int_{Y}Z_{1}^{\pi }\left( y\right) Z_{2}^{\pi }\left( y\right) \mu \left(
dy\right) \overset{law}{=}\int_{Y}Z_{1}^{\sigma }\left( y\right)
Z_{2}^{\sigma }\left( y\right) \mu \left( dy\right) \text{, for every }\left[
\pi \right] ,\left[ \sigma \right] \in \widehat{G}\text{.}
\label{identity in law}
\end{equation}

\bigskip

In the sequel, an identity such as (\ref{identity in law}) will be called a
\textbf{Watson's type relation}.

\section{Watson's type relations for Volterra processes}

\subsection{Setup and statement of the main results}

Throughout this section, $G$ stands for a finite group such that the
character $\chi _{\pi }\left( \cdot \right) $ is real-valued for every $%
\left[ \pi \right] \in \widehat{G}$. To simplify some technical points of
our discussion (in particular, to apply several crucial properties of
multiple Wiener-It\^{o} integrals) we will consider a two-dimensional,
real-valued Gaussian process $\left( Z_{1},Z_{2}\right) $ such that its
components are correlated Volterra processes. To define such objects, take a
measurable space $\left( T,\mathcal{T},\tau \right) $, where $\tau $ is
positive, $\sigma $-finite and non-atomic, and write $L_{%
\mathbb{R}
}^{2}\left( d\tau \right) $ to indicate the Hilbert space of real-valued,
square-integrable functions with respect to $\tau $. In what follows, we
will write
\begin{equation}
X=\left\{ X\left( f\right) :f\in L_{%
\mathbb{R}
}^{2}\left( d\tau \right) \right\}  \label{isonormal}
\end{equation}%
to indicate an \textbf{isonormal Gaussian process }(or a \textbf{Gaussian
measure}) on $L_{%
\mathbb{R}
}^{2}\left( d\tau \right) $. This means that $X$ is a centered Gaussian
family indexed by the elements of $L_{%
\mathbb{R}
}^{2}\left( d\tau \right) $, defined on some probability space $\left(
\Omega ,\mathcal{F},\mathbf{P}\right) $ and such that, for every $%
f_{1},f_{2}\in L_{%
\mathbb{R}
}^{2}\left( d\tau \right) $,%
\begin{equation*}
\mathbf{E}\left( X\left( f_{1}\right) X\left( f_{2}\right) \right)
=\int_{T}f_{1}\left( t\right) f_{2}\left( t\right) \tau \left( dt\right)
\text{.}
\end{equation*}

\bigskip

Now fix a measurable space $\left( Y,\mathcal{Y}\right) $. A two-dimensional
Gaussian process $\left\{ \left( Z_{1}\left( y\right) ,Z_{2}\left( y\right)
\right) :y\in Y\right\} $ is called a \textbf{correlated (generalized)
Volterra process},\textbf{\ }with respect to $X$ and with parameter $\rho
\in \left[ 0,1\right] $, if there exist two $\mathcal{Y}\otimes \mathcal{T}$
- measurable applications
\begin{equation*}
Y\times T\mapsto \mathbb{R}:\left( y,t\right) \mapsto \phi _{i}\left(
y,t\right) \text{, \ \ }i=1,2\text{,}
\end{equation*}%
such that: (a) for every $y\in Y$ the application $t\mapsto \phi _{i}\left(
y,t\right) $ is an element of $L_{%
\mathbb{R}
}^{2}\left( d\tau \right) $, (b) a.s. -- $\mathbb{P}$,%
\begin{equation}
Z_{i}\left( y\right) =X\left( \phi _{i}\left( y,\cdot \right) \right) \text{%
, \ \ }i=1,2\text{,}  \label{Volterra}
\end{equation}%
and (c) for every $y_{1},y_{2}\in Y$ and by using the notation introduced in
(\ref{COV}),%
\begin{eqnarray}
R_{1,1}\left( y_{1},y_{2}\right) &=&R_{2,2}\left( y_{1},y_{2}\right) \text{
\ and}  \label{correlated process} \\
R_{1,2}\left( y_{1},y_{2}\right) &=&R_{2,1}\left( y_{1},y_{2}\right) =\rho
R_{1,1}\left( y_{1},y_{2}\right) .  \notag
\end{eqnarray}

Note that, if $\rho =1$, then $Z_{1}\left( y\right) =Z_{2}\left( y\right) $
p.s.-$\mathbf{P}$, $\forall y\in Y$; moreover, the covariance structure of a
Gaussian process $\left( Z_{1},Z_{2}\right) $ of the type (\ref{Volterra})
may be rewritten as%
\begin{equation}
R_{i,j}\left( y_{1},y_{2}\right) =\mathbf{E}\left[ Z_{i}\left( y_{1}\right)
Z_{j}\left( y_{2}\right) \right] =\int_{T}\phi _{i}\left( y_{1},t\right)
\phi _{j}\left( y_{2},t\right) \tau \left( dt\right) \text{, \ \ }i,j=1,2;
\label{isoncov}
\end{equation}%
as a consequence, in view of (\ref{correlated process}) and (\ref{isoncov}),
and given an action $g\cdot y$ of $G$ on $Y$, $\left( Z_{1},Z_{2}\right) $
has a $G$-invariant law if, and only if, for $i$ equal to $1$ or $2$,%
\begin{equation}
\int_{T}\phi _{i}\left( g\cdot y_{1},t\right) \phi _{i}\left( g\cdot
y_{2},t\right) \tau \left( dt\right) =\int_{T}\phi _{i}\left( y_{1},t\right)
\phi _{i}\left( y_{2},t\right) \tau \left( dt\right) ,  \label{isonGINV}
\end{equation}%
for every $y_{1},y_{2}\in Y$ and every $g\in G$. In the sequel, to simplify
the notation, we will write
\begin{equation}
R_{1,1}\left( \cdot ,\cdot \right) =R_{2,2}\left( \cdot ,\cdot \right)
=R\left( \cdot ,\cdot \right) .  \label{simplification}
\end{equation}

\bigskip

We now fix an action $g\cdot y$ of $G$ on $Y$, as well as a $G$-invariant,
positive and $\sigma $-finite measure $\mu $ on $\left( Y,\mathcal{Y}\right)
$. For every real-valued $\Phi ,\Psi \in L^{2}\left( Y^{2},\mathcal{Y}%
^{2},d\mu \times d\mu \right) :=L^{2}\left( d\mu \times d\mu \right) $, we
define, for $y_{1},y_{2}\in Y$,

\begin{description}
\item[(i)] $\left[ \Phi \otimes _{\left( 1\right) }\Phi \right] \left(
y_{1},y_{2}\right) =\Phi \left( y_{1},y_{2}\right) $;

\item[(ii)] $\left[ \Phi \otimes _{\left( 2\right) }\Psi \right] \left(
y_{1},y_{2}\right) =\int_{Y}\Phi \left( y_{1},x\right) \Psi \left(
y_{2},x\right) \mu \left( dx\right) $;

\item[(iii)] $\forall p\geq 3$, $\left[ \Phi \otimes _{\left( p\right) }\Phi %
\right] \left( y_{1},y_{2}\right) =\left[ \left[ \Phi \otimes _{\left(
p-1\right) }\Phi \right] \otimes _{\left( 2\right) }\Phi \right] \left(
y_{1},y_{2}\right) ;$
\end{description}

\bigskip

Observe that, if $\Phi \in L^{2}\left( d\mu \times d\mu \right) $, then the
application $y\mapsto \left[ \Phi \otimes _{\left( p\right) }\Phi \right]
\left( y,y\right) $ is an element of $L^{1}\left( Y,\mathcal{Y},d\mu \right)
$ for every $p\geq 2$. Finally, for $\rho \in \left[ 0,1\right] $ as above,
we introduce the following set of real constants%
\begin{eqnarray}
K\left( 1,\rho \right) &=&2\rho  \notag \\
K\left( n,\rho \right) &=&2\sum_{j=0}^{\frac{n}{2}-1}\dbinom{n-1}{2j}\rho
^{2j}+2\sum_{j=0}^{\frac{n}{2}-1}\dbinom{n-1}{2j+1}\rho ^{2j+2}\text{, \ }n%
\text{ even,\ }n\geq 2\text{,}  \label{coeffs} \\
K\left( n,\rho \right) &=&2\sum_{j=0}^{\frac{n-1}{2}}\dbinom{n-1}{2j}\rho
^{2j+1}+2\sum_{j=0}^{\frac{n-3}{2}}\dbinom{n-1}{2j+1}\rho ^{2j+1}\text{, \ }n%
\text{ odd,\ }n\geq 3\text{.}  \notag
\end{eqnarray}

Note that $K\left( n,1\right) =2^{n}$ for every $n\geq 1$, $K\left( 2p,\rho
\right) >0$ for every $p\geq 1$, and, for $p\geq 0$, $K\left( 2p+1,\rho
\right) =0$ if, and only if, $\rho =0$ (since $\rho $ is real). In the next
result, under some additional integrability assumptions, we state necessary
and sufficient conditions to have that property (\ref{identity in law}) is
satisfied.

\bigskip

\begin{theorem}
Consider a finite group $G$ such that $\chi _{\pi }\left( \cdot \right) \in
\mathbb{R}$, for every $\left[ \pi \right] \in \widehat{G}$. Let the process
$\left( Z_{1},Z_{2}\right) $ be a correlated Volterra process of the type (%
\ref{Volterra}), for some correlation coefficient $\rho \in \left[ 0,1\right]
$, and assume $\left( Z_{1},Z_{2}\right) $ has a $G$-invariant law. Let also
$\mu \left( \cdot \right) $ be a $G$-invariant, positive measure satisfying
the assumptions of Proposition 7, and suppose moreover
\begin{equation}
\mathbf{E}\left( \int_{Y}Z_{1}\left( y\right) ^{2}\mu \left( dy\right)
\right) =\int_{Y}\int_{T}\phi _{1}\left( y,t\right) ^{2}\mu \left( dy\right)
\tau \left( dt\right) =\int_{Y}R\left( y,y\right) \mu \left( dy\right)
<+\infty \text{.}  \label{AssT8}
\end{equation}

Then,

\begin{enumerate}
\item the covariance functions $R$ and $R^{\pi \otimes \pi }$, defined
respectively according to (\ref{simplification}) and (\ref{reprCov}), for $%
\left[ \pi \right] \in \widehat{G}$, satisfy%
\begin{equation}
\int_{Y}\int_{Y}R\left( x,y\right) ^{2}\mu \left( dx\right) \mu \left(
dy\right) <+\infty \text{ \ \ and \ \ }\int_{Y}\int_{Y}R^{\pi \otimes \pi
}\left( x,y\right) ^{2}\mu \left( dx\right) \mu \left( dy\right) <+\infty ;
\label{sqintcov}
\end{equation}

\item the random variables
\begin{equation*}
\int_{Y}Z_{1}^{\pi }\left( y\right) Z_{2}^{\pi }\left( y\right) \mu \left(
dy\right) \text{, \ \ }\left[ \pi \right] \in \widehat{G}\text{,}
\end{equation*}%
are stochastically independent;

\item for every $\left[ \pi \right] \in \widehat{G}$, the process $\left(
Z_{1}^{\pi },Z_{2}^{\pi }\right) $ is a correlated Volterra process, with
parameter $\rho $;

\item the following three conditions are equivalent: (i) for every $\left[
\pi \right] ,\left[ \sigma \right] \in \widehat{G}$,%
\begin{equation}
\int_{Y}Z_{1}^{\pi }\left( y\right) Z_{2}^{\pi }\left( y\right) \mu \left(
dy\right) \overset{law}{=}\int_{Y}Z_{1}^{\sigma }\left( y\right)
Z_{2}^{\sigma }\left( y\right) \mu \left( dy\right) ,  \label{C I}
\end{equation}%
(ii) for every $\left[ \pi \right] \in \widehat{G}$ and every $n\geq 1$%
\begin{equation}
K\left( n,\rho \right) \int_{Y}\left[ R^{\pi \otimes \pi }\otimes _{\left(
n\right) }R^{\pi \otimes \pi }\right] \left( y,y\right) \mu \left( dy\right)
=\frac{K\left( n,\rho \right) }{\mid \widehat{G}\mid }\int_{Y}\left[
R\otimes _{\left( n\right) }R\right] \left( y,y\right) \mu \left( dy\right) ,
\label{C II}
\end{equation}%
(iii) for every $\left[ \pi \right] \in \widehat{G}$ and every $n\geq 1$%
\begin{equation}
K\left( n,\rho \right) \int_{Y}\left[ R^{\pi \otimes \pi }\otimes _{\left(
n\right) }R^{\pi \otimes \pi }\right] \left( y,y\right) \mu \left( dy\right)
=K\left( n,\rho \right) \int_{Y}\left[ R^{\sigma \otimes \sigma }\otimes
_{\left( n\right) }R^{\sigma \otimes \sigma }\right] \left( y,y\right) \mu
\left( dy\right) .  \label{C III}
\end{equation}
\end{enumerate}
\end{theorem}

\bigskip

\textbf{Remarks -- }(i) In view of (\ref{isoncov}), both formulae (\ref{C II}%
) and (\ref{C III}) can be immediately reformulated in terms of the kernels $%
\phi _{1}$ and $\phi _{2}$.

(ii) The role of the constants $K\left( n,\rho \right) $ in (\ref{C II}) and
(\ref{C III}) is immaterial for $\rho \neq 0$, or for $n$ even and $\rho \in %
\left[ 0,1\right] $.

\bigskip

Before proving Theorem 8, we state some interesting consequences of Theorem
8-4.

\bigskip

\begin{proposition}
Let $G=\left\{ e,g\right\} \simeq \mathbb{Z}/2\mathbb{Z}$, where $e$ stands
for the identity element. Keep the assumptions and the notation of Theorem
8, and suppose moreover that $\rho \neq 0$. Then, condition (\ref{C I}) is
verified if, and only if, for every $n\geq 1$%
\begin{equation*}
\int_{Y}\left[ R\otimes _{\left( n\right) }R\right] \left( y,g\cdot y\right)
\mu \left( dy\right) =0.
\end{equation*}
\end{proposition}

\begin{proof}
As already pointed out, in this case $\widehat{G}=\left\{ \left[ \pi _{u}%
\right] ,\left[ \pi _{a}\right] \right\} $, where $\left[ \pi _{u}\right] $
and $\left[ \pi _{a}\right] $ are the equivalence classes of the unity and
of the alternating representation. Moreover, due to Theorem 8 and the fact
that $K\left( n,\rho \right) \neq 0$, we know that (\ref{C I}) holds if, and
only if,
\begin{equation}
\int_{Y}\left[ R^{\pi _{u}\otimes \pi _{u}}\otimes _{\left( n\right) }R^{\pi
_{u}\otimes \pi _{u}}\right] \left( y,y\right) \mu \left( dy\right) =\frac{1%
}{2}\int_{Y}\left[ R\otimes _{\left( n\right) }R\right] \left( y,y\right)
\mu \left( dy\right) \text{,}  \label{CondZ/2Z}
\end{equation}%
for any $n\geq 1$, where
\begin{eqnarray*}
R^{\pi _{u}\otimes \pi _{u}}\left( y_{1},y_{2}\right) &=&\frac{1}{4}\left(
R\left( e\cdot y_{1},e\cdot y_{2}\right) +R\left( g\cdot y_{1},e\cdot
y_{2}\right) +R\left( e\cdot y_{1},g\cdot y_{2}\right) +R\left( g\cdot
y_{1},g\cdot y_{2}\right) \right) \\
&=&\frac{1}{2}\left( R\left( y_{1},y_{2}\right) +R\left( y_{1},g\cdot
y_{2}\right) \right) \text{,}
\end{eqnarray*}%
due to the $G$-invariance of the law of $\left( Z_{1},Z_{2}\right) $.
Finally, since $\mu $ is also $G$-invariant, one can easily prove that, for $%
n\geq 1$,%
\begin{equation*}
\int_{Y}\left[ R^{\pi _{u}\otimes \pi _{u}}\otimes _{\left( n\right) }R^{\pi
_{u}\otimes \pi _{u}}\right] \left( y,y\right) \mu \left( dy\right) =\frac{1%
}{2}\int_{Y}\left[ R\otimes _{\left( n\right) }R\right] \left( y,y\right)
\mu \left( dy\right) +\frac{1}{2}\int_{Y}\left[ R\otimes _{\left( n\right) }R%
\right] \left( y,g\cdot y\right) \mu \left( dy\right) ,
\end{equation*}%
thus yielding, via (\ref{CondZ/2Z}), the desired conclusion.
\end{proof}

\bigskip

\textbf{Remark -- }The process $\left( v_{1},v_{2}\right) $ defined in
formula (\ref{compensatedBr}) of the previous section can be represented as
as a correlated Volterra process, with parameter $\rho \in \left[ 0,1\right]
$ and covariance structure%
\begin{eqnarray}
\mathbf{E}\left[ v_{1}\left( s\right) v_{1}\left( t\right) \right] &=&%
\mathbf{E}\left[ v_{2}\left( s\right) v_{2}\left( t\right) \right]
=R_{v}\left( s,t\right) =s\wedge t-\frac{s+t}{2}+\frac{\left( s-t\right) ^{2}%
}{2}+\frac{1}{12}  \label{CovCompBr} \\
\mathbf{E}\left[ v_{1}\left( s\right) v_{2}\left( t\right) \right] &=&%
\mathbf{E}\left[ v_{2}\left( s\right) v_{1}\left( t\right) \right] =\rho
R_{v}\left( s,t\right) \text{,}  \notag
\end{eqnarray}%
where $s,t\in \left[ 0,1\right] $. Moreover, its law is $G$-invariant, where
$G=\mathbb{Z}/2\mathbb{Z}$. Since (\ref{WatRel}) holds, we deduce from
Corollary 9 that for every $n\geq 1$%
\begin{equation*}
\int_{0}^{1}\left[ R_{v}\otimes _{\left( n\right) }R_{v}\right] \left(
t,1-t\right) dt=0.
\end{equation*}

\bigskip

The next result, which is again a consequence of Theorem 8, is very useful
to deal with multiparameter processes.

\bigskip

\begin{proposition}
Fix $d\geq 2$. Let $\left( Y^{\left( i\right) },\mathcal{Y}^{\left( i\right)
},\mu ^{\left( i\right) }\right) $, $i=1,...,d$, be a collection of measure
spaces, with $\mu ^{\left( i\right) }$ positive and $\sigma $-finite, and
let $G^{\left( i\right) },...,G^{\left( d\right) }$ be finite groups with
real-valued characters, such that, for each $i=1,...,d$, an action $%
g_{i}\cdot y_{i}$ of $G_{i}$ on $Y_{i}$ is well defined. We note
\begin{eqnarray*}
Y &=&Y^{\left( 1\right) }\times \cdot \cdot \cdot \times Y^{\left( d\right) }%
\text{ \ , \ \ }G=G^{\left( 1\right) }\times \cdot \cdot \cdot \times
G^{\left( d\right) }\text{,} \\
\mu &=&\mu ^{\left( 1\right) }\times \cdot \cdot \cdot \times \mu ^{\left(
d\right) }
\end{eqnarray*}%
and we endow $G$ with the product group structure (see \cite[Section 3.2]%
{Serre}). Let also $\left\{ \left( Z_{1}\left( y\right) ,Z_{2}\left(
y\right) \right) :y\in Y\right\} $ be a correlated Volterra process with
parameter $\rho \in \left[ 0,1\right] $, such that, for every $y=\left(
y_{1},...,y_{d}\right) $ and $x=\left( x_{1},...,x_{d}\right) $ in $Y$%
\begin{equation*}
\mathbf{E}\left[ Z_{1}\left( x\right) Z_{1}\left( y\right) \right] =R\left(
x,y\right) =\prod_{i=1}^{d}R_{i}\left( x_{i},y_{i}\right) \text{ \ \ and \ \
}\mathbf{E}\left[ Z_{1}\left( x\right) Z_{2}\left( y\right) \right] =\rho
R\left( x,y\right) \text{,}
\end{equation*}%
where for each $i$, $R_{i}$ is a $G^{\left( i\right) }$-invariant covariance
function such that%
\begin{equation*}
\int_{Y^{\left( i\right) }}R_{i}\left( y,y\right) \mu ^{\left( i\right)
}\left( dy\right) <+\infty .
\end{equation*}%
Then,

\begin{enumerate}
\item the application
\begin{equation*}
\left( g,y\right) \mapsto g\cdot y:\left(
g_{1},...,g_{d};y_{1},...,y_{d}\right) \mapsto \left( g_{1}\cdot
y_{1},...,g_{d}\cdot y_{d}\right)
\end{equation*}%
is an action of $G$ on $Y$;

\item the process $\left( Z_{1},Z_{2}\right) $ has a $G$-invariant law;

\item for every $\left[ \pi _{i}\right] ,\left[ \sigma _{i}\right] \in
\widehat{G^{\left( i\right) }}$, $i=1,...,d$, every $x,y\in Y$%
\begin{equation*}
R^{\theta }\left( x,y\right) =\left\{
\begin{array}{lll}
0 &  & \text{if there exists }i\text{ such that }\left[ \pi _{i}\right] \neq %
\left[ \sigma _{i}\right] \\
\prod_{i=1}^{d}R^{\pi _{i}\otimes \pi _{i}}\left( x_{i},y_{i}\right) &  &
\text{otherwise,}%
\end{array}%
\right.
\end{equation*}%
where $\theta =\left[ \pi _{i}\otimes \cdot \cdot \cdot \otimes \pi _{d}%
\right] \otimes \left[ \sigma _{i}\otimes \cdot \cdot \cdot \otimes \sigma
_{d}\right] $ is a generic element of $\widehat{G^{2}}$;

\item if, for each $i=1,...,d$, the function $R_{i}$ satisfies either one of
conditions (\ref{C II}) and (\ref{C III}), then $\left( Z_{1},Z_{2}\right) $
verifies Watson's relation (\ref{C I}) for every $\left[ \pi \right] ,\left[
\sigma \right] \in \widehat{G}$.
\end{enumerate}
\end{proposition}

\begin{proof}
Points 1. and 2.\ are trivial. Point 3 is a consequence of the $G^{\left(
i\right) }$ invariance of each $R_{i}$, as well as of Proposition 6-3. To
prove point 4., suppose that each $R_{i}$ verifies (\ref{C III}), and also $%
\rho \neq 0$. Then, $K\left( n,\rho \right) \neq 0$ for each $n$, and%
\begin{eqnarray*}
\frac{1}{\mid \widehat{G}\mid }\int_{Y}\left[ R\otimes _{\left( n\right) }R%
\right] \left( y,y\right) \mu \left( dy\right) &=&\frac{1}{\mid \widehat{G}%
\mid }\prod_{i=1}^{d}\int_{Y^{\left( i\right) }}\left[ R_{i}\otimes _{\left(
n\right) }R_{i}\right] \left( y_{i},y_{i}\right) \mu ^{\left( i\right)
}\left( dy_{i}\right) \\
&=&\prod_{i=1}^{d}\int_{Y^{\left( i\right) }}\left[ R_{i}^{\left[ \pi _{i}%
\right] \otimes \left[ \pi _{i}\right] }\otimes _{\left( n\right) }R_{i}^{%
\left[ \pi _{i}\right] \otimes \left[ \pi _{i}\right] }\right] \left(
y_{i},y_{i}\right) \mu ^{\left( i\right) }\left( dy_{i}\right) \text{,}
\end{eqnarray*}%
for every $\left[ \pi _{i}\right] \in \widehat{G^{\left( i\right) }}$, $%
i=1,...,d$, since $\mid \widehat{G}\mid =\prod_{i=1,...,d}\mid \widehat{%
G^{\left( i\right) }}\mid $. To conclude, just observe that, thanks to point
3.,%
\begin{equation*}
\prod_{i=1}^{d}\int_{Y^{\left( i\right) }}\left[ R_{i}^{\left[ \pi _{i}%
\right] \otimes \left[ \pi _{i}\right] }\otimes _{\left( n\right) }R_{i}^{%
\left[ \pi _{i}\right] \otimes \left[ \pi _{i}\right] }\right] \left(
y_{i},y_{i}\right) \mu ^{\left( i\right) }\left( dy_{i}\right) =\int_{Y}%
\left[ R^{\eta \otimes \eta }\otimes _{\left( n\right) }R^{\eta \otimes \eta
}\right] \left( y,y\right) \mu \left( dy\right) \text{,}
\end{equation*}%
where $\eta =\pi _{i}\otimes \cdot \cdot \cdot \otimes \pi _{d}$ (to deal
with the case $\rho =0$, just perform the same argument for even cumulants).
\end{proof}

\bigskip

\textbf{Example }(A quadruplication identity). Let $\mathbb{B}_{0}=\left\{
\mathbb{B}_{0,1}\left( t_{1},t_{2}\right) ,\mathbb{B}_{0,2}\left(
t_{1},t_{2}\right) :\left( t_{1},t_{2}\right) \in \left[ 0,1\right]
^{2}\right\} $ be a correlated tied-down Brownian sheet, that is, $\mathbb{B}%
_{0}$ is a two-dimensional Gaussian process such that
\begin{eqnarray*}
\mathbf{E}\left[ \mathbb{B}_{0,1}\left( t_{1},t_{2}\right) ,\mathbb{B}%
_{0,1}\left( s_{1},s_{2}\right) \right] &=&\mathbf{E}\left[ \mathbb{B}%
_{0,2}\left( t_{1},t_{2}\right) ,\mathbb{B}_{0,2}\left( s_{1},s_{2}\right) %
\right] \\
&=&\left( t_{1}\wedge s_{1}-s_{1}t_{1}\right) \left( t_{2}\wedge
s_{2}-s_{2}t_{2}\right) \text{,} \\
\mathbf{E}\left[ \mathbb{B}_{0,1}\left( t_{1},t_{2}\right) ,\mathbb{B}%
_{0,2}\left( s_{1},s_{2}\right) \right] &=&\mathbf{E}\left[ \mathbb{B}%
_{0,2}\left( t_{1},t_{2}\right) ,\mathbb{B}_{0,1}\left( s_{1},s_{2}\right) %
\right] \\
&=&\rho \times \left( t_{1}\wedge s_{1}-s_{1}t_{1}\right) \left( t_{2}\wedge
s_{2}-s_{2}t_{2}\right) \text{,}
\end{eqnarray*}%
where $\left( s_{1},s_{2}\right) ,\left( t_{1},t_{2}\right) \in \left[ 0,1%
\right] ^{2}$ and $\rho \in \left[ 0,1\right] $. Note that $\mathbb{B}_{0}$
can be represented as a Volterra process. Moreover, $\mathbb{B}_{0}$ has the
law of a correlated Brownian sheet $\mathbb{W}$ (with the same parameter%
\footnote{%
That is,
\begin{equation*}
\mathbb{W}=\left\{ \mathbb{W}_{1}\left( t_{1},t_{2}\right) ,\mathbb{W}%
_{2}\left( t_{1},t_{2}\right) :\left( t_{1},t_{2}\right) \in \left[ 0,1%
\right] ^{2}\right\}
\end{equation*}%
is a two-dimensional Gaussian process such that
\begin{eqnarray*}
\mathbf{E}\left[ \mathbb{W}_{1}\left( t_{1},t_{2}\right) ,\mathbb{W}%
_{1}\left( s_{1},s_{2}\right) \right] &=&\mathbf{E}\left[ \mathbb{W}%
_{2}\left( t_{1},t_{2}\right) ,\mathbb{W}_{2}\left( s_{1},s_{2}\right) %
\right] =\left( t_{1}\wedge s_{1}\right) \left( t_{2}\wedge s_{2}\right)
\text{,} \\
\mathbf{E}\left[ \mathbb{W}_{1}\left( t_{1},t_{2}\right) ,\mathbb{W}%
_{2}\left( s_{1},s_{2}\right) \right] &=&\mathbf{E}\left[ \mathbb{W}%
_{2}\left( t_{1},t_{2}\right) ,\mathbb{W}_{1}\left( s_{1},s_{2}\right) %
\right] =\rho \times \left( t_{1}\wedge s_{1}\right) \left( t_{2}\wedge
s_{2}\right) \text{.}
\end{eqnarray*}%
}), conditioned to vanish on the edges of the square $\left[ 0,1\right] ^{2}$%
. Now define, for $i=1,2$, the compensated processes%
\begin{equation*}
\mathbb{U}_{i}\left( t_{1},t_{2}\right) =\mathbb{B}_{0,i}\left(
t_{1},t_{2}\right) -\int_{0}^{1}\mathbb{B}_{0,i}\left( t_{1},u_{2}\right)
du_{2}-\int_{0}^{1}\mathbb{B}_{0,i}\left( u_{1},t_{2}\right) du_{1}+\int_{%
\left[ 0,1\right] ^{2}}\mathbb{B}_{0,i}\left( u_{1},u_{2}\right) du_{1}du_{2}%
\text{,}
\end{equation*}%
where $\left( t_{1},t_{2}\right) \in \left[ 0,1\right] ^{2}$. We claim that
the following identity in law holds%
\begin{equation}
\int_{\left[ 0,1\right] ^{2}}\mathbb{U}_{1}\left( t_{1},t_{2}\right) \mathbb{%
U}_{2}\left( t_{1},t_{2}\right) dt_{1}dt_{2}\overset{law}{=}\frac{1}{16}%
\sum_{i=1}^{4}\int_{\left[ 0,1\right] ^{2}}\mathbb{B}_{0,1}^{\left( i\right)
}\left( t_{1},t_{2}\right) \mathbb{B}_{0,2}^{\left( i\right) }\left(
t_{1},t_{2}\right) \left( t_{1},t_{2}\right) ^{2}dt_{1}dt_{2}\text{,}
\label{Quad}
\end{equation}%
where $\mathbb{B}_{0}^{\left( i\right) }=\left( \mathbb{B}_{0,1}^{\left(
i\right) },\mathbb{B}_{0,2}^{\left( i\right) }\right) $, $i=1,...,4$, are
four independent copies of $\mathbb{B}_{0}$. As a matter of fact, standard
calculations show that%
\begin{eqnarray*}
\mathbf{E}\left[ \mathbb{U}_{1}\left( t_{1},t_{2}\right) \mathbb{U}%
_{1}\left( s_{1},s_{2}\right) \right] &=&\mathbf{E}\left[ \mathbb{U}%
_{2}\left( t_{1},t_{2}\right) \mathbb{U}_{2}\left( s_{1},s_{2}\right) \right]
=R_{v}\left( s_{1},t_{1}\right) R_{v}\left( s_{2},t_{2}\right) \text{,} \\
\mathbf{E}\left[ \mathbb{U}_{1}\left( t_{1},t_{2}\right) \mathbb{U}%
_{2}\left( s_{1},s_{2}\right) \right] &=&\mathbf{E}\left[ \mathbb{U}%
_{2}\left( t_{1},t_{2}\right) \mathbb{U}_{1}\left( s_{1},s_{2}\right) \right]
=\rho \times R_{v}\left( s_{1},t_{1}\right) R_{v}\left( s_{2},t_{2}\right)
\text{,}
\end{eqnarray*}%
where $R_{v}$ is defined as in (\ref{CovCompBr}). Since $R_{v}$ is invariant
with respect to the action of $\left\{ e,g\right\} \simeq \left( \mathbb{Z}/2%
\mathbb{Z}\right) $ on $\left[ 0,1\right] $ given by $e\cdot t=t$ and $%
g\cdot t=1-t$, Proposition 10-2 entails that the law of the vector $\left(
\mathbb{U}_{1},\mathbb{U}_{2}\right) $ is invariant with respect to the
action of the product group $G=\left\{ e,g\right\} \times \left\{
e,g\right\} \simeq \left( \mathbb{Z}/2\mathbb{Z}\right) \times \left(
\mathbb{Z}/2\mathbb{Z}\right) $ on $\left[ 0,1\right] ^{2}$ defined as%
\begin{eqnarray*}
\left( e,e\right) \cdot \left( t_{1},t_{2}\right) &=&\left(
t_{1},t_{2}\right) \text{, }\left( e,g\right) \cdot \left(
t_{1},t_{2}\right) =\left( t_{1},1-t_{2}\right) \\
\left( g,e\right) \cdot \left( 1-t_{1},t_{2}\right) &=&\left(
t_{1},t_{2}\right) \text{, }\left( g,g\right) \cdot \left(
t_{1},t_{2}\right) =\left( 1-t_{1},1-t_{2}\right) .
\end{eqnarray*}%
Now recall that $\widehat{G}=\left\{ \left[ \pi _{u}\right] \otimes \left[
\pi _{u}\right] ,\left[ \pi _{a}\right] \otimes \left[ \pi _{u}\right] ,%
\left[ \pi _{u}\right] \otimes \left[ \pi _{a}\right] ,\left[ \pi _{a}\right]
\otimes \left[ \pi _{a}\right] \right\} $, where $\pi _{u}$ and $\pi _{a}$
are the unity and \textit{alternating} representation of $\mathbb{Z}/2%
\mathbb{Z}$. According to Proposition 7-4 (since Lebesgue measure on $\left[
0,1\right] ^{2}$ is also $G$-invariant) and Proposition 10-4, for every $%
\lambda \in \mathbb{R}$,
\begin{eqnarray*}
&&\mathbf{E}\left[ \exp \left( \mathtt{i}\lambda \int_{\left[ 0,1\right]
^{2}}\mathbb{U}_{1}\left( t_{1},t_{2}\right) \mathbb{U}_{2}\left(
t_{1},t_{2}\right) dt_{1}dt_{2}\right) \right] \\
&=&\mathbf{E}\left[ \exp \left( \mathtt{i}\lambda \int_{\left[ 0,1\right]
^{2}}\mathbb{U}_{1}^{\pi _{a}\otimes \pi _{a}}\left( t_{1},t_{2}\right)
\mathbb{U}_{2}^{\pi _{a}\otimes \pi _{a}}\left( t_{1},t_{2}\right)
dt_{1}dt_{2}\right) \right] ^{4}.
\end{eqnarray*}%
To conclude, we use Proposition 7-2 to show that
\begin{eqnarray*}
\mathbf{E}\left[ \mathbb{U}_{1}^{\pi _{a}\otimes \pi _{a}}\left(
t_{1},t_{2}\right) \mathbb{U}_{1}^{\pi _{a}\otimes \pi _{a}}\left(
s_{1},s_{2}\right) \right] &=&\mathbf{E}\left[ \mathbb{U}_{2}^{\pi
_{a}\otimes \pi _{a}}\left( t_{1},t_{2}\right) \mathbb{U}_{2}^{\pi
_{a}\otimes \pi _{a}}\left( s_{1},s_{2}\right) \right] \\
&=&R_{v}^{\pi _{a}\otimes \pi _{a}}\left( s_{1},t_{1}\right) R_{v}^{\pi
_{a}\otimes \pi _{a}}\left( s_{2},t_{2}\right) \text{,} \\
\mathbf{E}\left[ \mathbb{U}_{1}^{\pi _{a}\otimes \pi _{a}}\left(
t_{1},t_{2}\right) \mathbb{U}_{2}^{\pi _{a}\otimes \pi _{a}}\left(
s_{1},s_{2}\right) \right] &=&\mathbf{E}\left[ \mathbb{U}_{2}^{\pi
_{a}\otimes \pi _{a}}\left( t_{1},t_{2}\right) \mathbb{U}_{1}^{\pi
_{a}\otimes \pi _{a}}\left( s_{1},s_{2}\right) \right] \\
&=&\rho \times R_{v}^{\pi _{a}\otimes \pi _{a}}\left( s_{1},t_{1}\right)
R_{v}^{\pi _{a}\otimes \pi _{a}}\left( s_{2},t_{2}\right) \text{,}
\end{eqnarray*}%
thus implying that
\begin{eqnarray*}
&&\left\{ \mathbb{U}_{1}\left( t_{1},t_{2}\right) ,\mathbb{U}_{2}\left(
t_{1},t_{2}\right) :\left( t_{1},t_{2}\right) \in \left[ 0,1/2\right]
^{2}\right\} \\
&&\overset{law}{=}\left\{ 4^{-1}\mathbb{B}_{0,1}\left( 2t_{1},2t_{2}\right)
,4^{-1}\mathbb{B}_{0,2}\left( 2t_{1},2t_{2}\right) :\left(
t_{1},t_{2}\right) \in \left[ 0,1/2\right] ^{2}\right\} \text{,}
\end{eqnarray*}%
and therefore
\begin{equation*}
\int_{\left[ 0,1\right] ^{2}}\mathbb{U}_{1}^{\pi _{a}\otimes \pi _{a}}\left(
t_{1},t_{2}\right) \mathbb{U}_{2}^{\pi _{a}\otimes \pi _{a}}\left(
t_{1},t_{2}\right) dt_{1}dt_{2}\overset{law}{=}\frac{1}{4}\int_{\left[ 0,1/2%
\right] ^{2}}\mathbb{B}_{0,1}\left( 2t_{1},2t_{2}\right) \mathbb{B}%
_{0,2}\left( 2t_{1},2t_{2}\right) dt_{1}dt_{2}\text{,}
\end{equation*}%
so that (\ref{Quad}) is obtained by a standard change of variables on the
right hand side of the previous expression. The reader is referred to \cite%
{PeYor2005} for other two-parameters generalizations of Watson identity.

\subsection{Proof of Theorem 8}

\textbf{(1.) }Since $G$ is finite, to prove both inequalities in formula (%
\ref{sqintcov}) it is sufficient to show that, for every $g,h\in G$,%
\begin{equation*}
\int_{Y}\int_{Y}R\left( h\cdot y,g\cdot z\right) ^{2}\mu \left( dz\right)
\mu \left( dy\right) <+\infty .
\end{equation*}%
But, since $\mu $ is $G$-invariant, and taking into account (\ref{isoncov}),%
\begin{eqnarray*}
\int_{Y}\int_{Y}R\left( h\cdot y,g\cdot z\right) ^{2}\mu \left( dz\right)
\mu \left( dy\right) &=&\int_{Y}\int_{Y}R\left( y,z\right) ^{2}\mu \left(
dz\right) \mu \left( dy\right) \\
&=&\int_{Y}\int_{Y}\left( \int_{T}\phi _{1}\left( z,t\right) \phi _{1}\left(
y,t\right) \tau \left( dt\right) \right) ^{2}\mu \left( dz\right) \mu \left(
dy\right) \\
&<&+\infty ,
\end{eqnarray*}%
due to (\ref{AssT8}), as well as to an application of the Cauchy-Schwarz
inequality.

\noindent%
\textbf{(2.) }This is a direct consequence of Proposition 7-1..

\noindent%
\textbf{(3.) }By additivity of Gaussian measures, for every $y\in Y$, $i=1,2$
and $\left[ \pi \right] \in \widehat{G}$,
\begin{equation*}
Z_{i}^{\pi }\left( y\right) =\frac{1}{\left\vert G\right\vert }\sum_{g\in
G}Z\left( g\cdot x\right) \chi _{\pi }\left( g^{-1}\right) =\frac{1}{%
\left\vert G\right\vert }\sum_{g\in G}X\left( \phi _{i}\left( g\cdot y,\cdot
\right) \right) \chi _{\pi }\left( g^{-1}\right) =X\left( \phi _{i}^{\left(
\pi \right) }\left( y,\cdot \right) \right)
\end{equation*}%
where
\begin{equation}
\phi _{i}^{\left( \pi \right) }\left( y,t\right) :=\frac{1}{\left\vert
G\right\vert }\sum_{g\in G}\phi _{i}\left( g\cdot y,t\right) \chi _{\pi
}\left( g^{-1}\right) \text{, \ \ }\left( y,t\right) \in Y\times T\text{.}
\label{kernelR}
\end{equation}%
(note that $\phi _{i}^{\left( \pi \right) }\in L^{2}\left( d\mu \times d\tau
\right) ).$ Moreover, for any $y_{1},y_{2}\in Y,$%
\begin{eqnarray*}
\mathbf{E}\left[ Z_{1}^{\pi }\left( y_{1}\right) Z_{1}^{\pi }\left(
y_{2}\right) \right] &=&\mathbf{E}\left[ Z_{2}^{\pi }\left( y_{1}\right)
Z_{2}^{\pi }\left( y_{2}\right) \right] =R^{\pi \otimes \pi }\left(
y_{1},y_{2}\right) \text{, \ \ and} \\
\mathbf{E}\left[ Z_{1}^{\pi }\left( y_{1}\right) Z_{2}^{\pi }\left(
y_{2}\right) \right] &=&\rho R^{\pi \otimes \pi }\left( y_{1},y_{2}\right) ,
\end{eqnarray*}%
due to formula (\ref{reprCov}), thus yielding the desired result.

\noindent%
\textbf{(4.) }Fix $\left[ \pi \right] \in \widehat{G}$. Since $Z_{1}^{\pi }$
and $Z_{2}^{\pi }$ are Volterra processes with respect to the Gaussian
measure $X$, we may apply a standard version of the multiplication formula
for Wiener-It\^{o} integrals (see e.g. \cite[p. 211]{DMM}) to obtain
\begin{equation*}
\int_{Y}Z_{1}^{\pi }\left( y\right) Z_{2}^{\pi }\left( y\right) \mu \left(
dy\right) =\int_{Y}\rho R^{\pi \otimes \pi }\left( y,y\right) \mu \left(
dy\right) +I_{2}^{X}\left( \Phi ^{\left( \pi \right) }\right) \text{,}
\end{equation*}%
where $I_{2}^{X}$ stands for a double Wiener-It\^{o} integral with respect
to $X$ (see again \cite{DMM}), and $\Phi ^{\left( \pi \right) }$ is the
symmetrized kernel%
\begin{equation*}
\Phi ^{\left( \pi \right) }\left( t_{1},t_{2}\right) =\frac{1}{2}\int_{Y}%
\left[ \phi _{1}^{\left( \pi \right) }\left( y,t_{1}\right) \phi
_{2}^{\left( \pi \right) }\left( y,t_{2}\right) +\phi _{1}^{\left( \pi
\right) }\left( y,t_{2}\right) \phi _{2}^{\left( \pi \right) }\left(
y,t_{1}\right) \right] \mu \left( dy\right) ,
\end{equation*}%
where $\phi _{i}^{\left( \pi \right) }$, $i=1,2$, is defined as in (\ref%
{kernelR}). On the other hand,
\begin{eqnarray*}
\int_{Y}Z_{1}\left( y\right) Z_{2}\left( y\right) \mu \left( dy\right)
&=&\int_{Y}\rho R\left( y,y\right) \mu \left( dy\right) +I_{2}^{X}\left(
\Phi \right) \text{, \ \ where } \\
\Phi \left( t_{1},t_{2}\right) &=&\frac{1}{2}\int_{Y}\left[ \phi _{1}\left(
y,t_{1}\right) \phi _{2}\left( y,t_{2}\right) +\phi _{1}\left(
y,t_{2}\right) \phi _{2}\left( y,t_{1}\right) \right] \mu \left( dy\right) .
\end{eqnarray*}%
Now, it is well known that the law of a double Wiener-It\^{o} integral is
determined by its cumulants (see \cite{Slud}). We therefore note $\kappa
_{n}\left( J\right) $, $n\geq 1$, the $n$th cumulant of a given random
variable $J$, and use a version of the \textit{diagram formulae} for
cumulants of multiple stochastic integrals (as presented, for instance, in
\cite{Surg}, \cite[Proposition 9 and Corollary 1]{Rota Wall} or \cite[%
Section 2]{Fox Taqqu}) to obtain that, for every $n\geq 2$, there exists a
universal combinatorial coefficient $c_{n}>0$ such that, for any $\left[ \pi %
\right] \in \widehat{G}$,%
\begin{equation*}
\kappa _{n}\left( I_{2}^{X}\left( \Phi ^{\left( \pi \right) }\right) \right)
=c_{n}\int_{T^{n}}\text{ }\Phi ^{\left( \pi \right) }\left(
t_{1},t_{2}\right) \Phi ^{\left( \pi \right) }\left( t_{2},t_{3}\right)
\cdot \cdot \cdot \Phi ^{\left( \pi \right) }\left( t_{n},t_{1}\right) \tau
\left( dt_{1}\right) \cdot \cdot \cdot \tau \left( dt_{n}\right)
\end{equation*}%
and also%
\begin{equation*}
\kappa _{n}\left( I_{2}^{X}\left( \Phi \right) \right) =c_{n}\int_{T^{n}}%
\text{ }\Phi \left( t_{1},t_{2}\right) \Phi \left( t_{2},t_{3}\right) \cdot
\cdot \cdot \Phi \left( t_{n},t_{1}\right) \tau \left( dt_{1}\right) \cdot
\cdot \cdot \tau \left( dt_{n}\right) .
\end{equation*}%
By using the relations%
\begin{eqnarray*}
\int_{T}\phi _{1}\left( y_{1},t\right) \phi _{1}\left( y_{2},t\right) \tau
\left( dt\right) &=&\int_{T}\phi _{2}\left( y_{1},t\right) \phi _{2}\left(
y_{2},t\right) \tau \left( dt\right) =R\left( y_{1},y_{2}\right) \\
\int_{T}\phi _{1}\left( y_{1},t\right) \phi _{2}\left( y_{2},t\right) \tau
\left( dt\right) &=&\int_{T}\phi _{2}\left( y_{1},t\right) \phi _{1}\left(
y_{2},t\right) \tau \left( dt\right) =\rho R\left( y_{1},y_{2}\right) \text{%
, }
\end{eqnarray*}%
as well as a combinatorial argument, we finally obtain, for $n\geq 2$,%
\begin{eqnarray*}
\kappa _{n}\left( I_{2}^{X}\left( \Phi ^{\left( \pi \right) }\right) \right)
&=&\frac{c_{n}\times K\left( n,\rho \right) }{2^{n}}\int_{Y}\left[ R^{\pi
\otimes \pi }\otimes _{\left( n\right) }R^{\pi \otimes \pi }\right] \left(
y,y\right) \mu \left( dy\right) \\
\kappa _{n}\left( I_{2}^{X}\left( \Phi \right) \right) &=&\frac{c_{n}\times
K\left( n,\rho \right) }{2^{n}}\int_{Y}\left[ R\otimes _{\left( n\right) }R%
\right] \left( y,y\right) \mu \left( dy\right) .
\end{eqnarray*}%
To conclude, use independence to write, for $n\geq 1$,%
\begin{equation*}
\kappa _{n}\left( \int_{Y}Z_{1}\left( y\right) Z_{2}\left( y\right) \mu
\left( dy\right) \right) =\sum_{\left[ \pi \right] \in \widehat{G}}\kappa
_{n}\left( \int_{Y}Z_{1}^{\pi }\left( y\right) Z_{2}^{\pi }\left( y\right)
\mu \left( dy\right) \right) \text{,}
\end{equation*}
and observe that, thanks to the translation invariance property of cumulants
(see e.g. \cite[Corollary 4.1]{RotaShen}), for any $n\geq 2$,%
\begin{eqnarray*}
\kappa _{n}\left( \int_{Y}Z_{1}\left( y\right) Z_{2}\left( y\right) \mu
\left( dy\right) \right) &=&\kappa _{n}\left( I_{2}^{X}\left( \Phi \right)
\right) \\
\kappa _{n}\left( \int_{Y}Z_{1}^{\pi }\left( y\right) Z_{2}^{\pi }\left(
y\right) \mu \left( dy\right) \right) &=&\kappa _{n}\left( I_{2}^{X}\left(
\Phi ^{\left( \pi \right) }\right) \right) \text{, \ \ }\left[ \pi \right]
\in \widehat{G}\text{, }
\end{eqnarray*}%
and consequently
\begin{equation*}
\kappa _{n}\left( I_{2}^{X}\left( \Phi \right) \right) =\sum_{\left[ \pi %
\right] \in \widehat{G}}\kappa _{n}\left( I_{2}^{X}\left( \Phi ^{\left( \pi
\right) }\right) \right) \text{.}
\end{equation*}%
The proof is completed by standard arguments. $\blacksquare $

\bigskip

\section{Refinements, further applications and examples}

\subsection{Connections with Karhunen-Lo\`{e}ve expansions}

In this paragraph, we elucidate some of the connections between
our decomposition of stochastic processes, and and
Karhunen-Lo\`{e}ve (KL) expansions of Gaussian processes indexed
by the elements of a measurable space $\left( T,\mathcal{T}\right)
$ (for fundamental facts about KL
expansions, see e.g. \cite{Adler}, \cite[Chapter 5]{SW}, as well as \cite%
{DPY} and the references therein). In what follows, $G$ is a topological
compact group, acting on $T$ through the application $\left( g,t\right)
\mapsto g\cdot t$, $t\in T$. We write $m\left( dt\right) $ to indicate a $G$
invariant measure on $\left( T,\mathcal{T}\right) $.

\bigskip

We also consider a \textit{positive definite} kernel $R\left( s,t\right) $, $%
s,t\in T$, such that $R$ is the covariance function of a centered
Gaussian stochastic process $\mathbb{X}=\left\{ \mathbb{X}\left(
t\right) :t\in
T\right\} $, defined on some probability space $\left( \Omega ,\mathcal{F},%
\mathbf{P}\right) $, and such that, for every $\omega \in \Omega $, the
function $t\mapsto \mathbb{X}\left( t\right) $ is in $L^{2}\left(
T,dm\right) $. We note $\lambda _{1}>\lambda _{2}>...>0$ the sequence of the
eigenvalues of $R$ (with respect to $m\left( \cdot \right) $), whereas $%
E_{1},E_{2},...$ indicate the associated eigenspaces. For every $j\geq 1$, $%
n_{j}$ is the (finite) dimension of $E_{j}$. The next assumption will be in
order throughout the paragraph

\bigskip

\textbf{Assumption C} -- For every $j$, note $\left(
f_{j,1},...,f_{j,n_{j}}\right) $ an orthonormal basis of $E_{j}$ (in the
sense of $L^{2}\left( T,dm\right) $). (C-i) The process $\mathbb{X}$ admits
the following KL expansion: there exists an array of i.i.d. $N\left(
0,1\right) $ random variables $\left\{ \xi _{j,l}:j\geq
1,l=1,...,n_{j}\right\} $ such that, as $N\rightarrow +\infty $ and p.s. - $%
\mathbf{P}$, the process%
\begin{equation*}
\mathbb{X}_{N}\left( t\right) =\sum_{j=1}^{N}\sqrt{\lambda _{j}}\left\{ \xi
_{j,1}\times f_{j,1}\left( t\right) +\xi _{j,2}\times f_{j,2}\left( t\right)
+\cdot \cdot \cdot +\xi _{j,n_{j}}\times f_{j,n_{j}}\left( t\right) \right\}
\text{, \ \ }t\in T\text{,}
\end{equation*}%
converges to $\mathbb{X}$ in $L^{2}\left( T,dm\right) $. (C-ii) The
processes $\mathbb{X}$, $\mathbb{X}_{N}$ ($N\geq 1$) and $f_{j,l}$ ($j\geq
1,l=1,...,n_{j}$) satisfy Assumption A of section 3.1, with $Y=T$.

\bigskip

The reader is referred once again to \cite{Adler} or \cite{SW} for (rather
general) sufficient conditions, ensuring the validity of Assumption (C-i) in
the case $\left[ 0,1\right] ^{d}$, $d\geq 1$. Assumption (C-ii) is redundant
for $G$ finite.

\bigskip

According to (\ref{projmarg2}), for every $\left[ \pi \right] \in \widehat{G}
$ we define%
\begin{eqnarray}
f^{\pi }\left( t\right) &=&d_{\pi }\int_{G}\chi _{\pi }\left( g\right)
f\left( g^{-1}\cdot t\right) dg\text{, \ \ }f\in L^{2}\left( T,dm\right)
\text{,}  \label{fpi} \\
\mathbb{X}_{N}^{\pi }\left( t\right) &=&d_{\pi }\int_{G}\chi _{\pi }\left(
g\right) \mathbb{X}_{N}\left( g^{-1}\cdot t\right) dg  \notag \\
&=&\sum_{j=1}^{N}\sqrt{\lambda _{j}}\left\{ \xi _{j,1}\times f_{j,1}^{\pi
}\left( t\right) +\xi _{j,2}\times f_{j,2}^{\pi }\left( t\right) +\cdot
\cdot \cdot +\xi _{j,n_{j}}\times f_{j,n_{j}}^{\pi }\left( t\right) \right\}
,  \notag \\
\mathbb{X}^{\pi }\left( t\right) &=&d_{\pi }\int_{G}\chi _{\pi }\left(
g\right) \mathbb{X}^{\pi }\left( g^{-1}\cdot t\right) dg.  \notag
\end{eqnarray}

Note that, according to Proposition 7-3, a.s.-$\mathbf{P}$,
\begin{equation}
\mathbb{X}\left( t\right) =\sum_{\left[ \pi \right] \in \widehat{G}}\mathbb{X%
}^{\pi }\left( t\right) \text{, \ \ }  \label{trajdec}
\end{equation}%
with convergence in $L^{2}\left( T,dm\right) $. The following Proposition
explains some remarkable relation between (\ref{trajdec}) and KL expansions,
in the case of Gaussian processes with a $G$ invariant law..

\bigskip

\begin{proposition}
Let the notation and assumptions of this paragraph prevail. Then,

\begin{enumerate}
\item for each $\left[ \pi \right] \in \widehat{G}$, $\mathbb{X}_{N}^{\pi
}\left( t\right) \rightarrow \mathbb{X}^{\pi }\left( t\right) $, as $%
N\rightarrow +\infty $, a.s.-$\mathbf{P}$ in $L^{2}\left( T,dm\right) ;$

\item suppose $\mathbb{X}$ has a $G$-invariant law; then, for each $j\geq 1$%
, the application
\begin{equation*}
g\mapsto \left\{ f\left( t\right) \mapsto f\left( g^{-1}\cdot t\right) :f\in
E_{j}\right\}
\end{equation*}%
is a finite dimensional representation of $G$;

\item for $j\geq 1$, write
\begin{equation*}
E_{j}=E_{j}^{1}\oplus \cdot \cdot \cdot \oplus E_{j}^{h_{j}}\text{,}
\end{equation*}%
with $1\leq h_{j}\leq n_{j}$, to indicate the \textbf{canonical decomposition%
} of $E_{j}$, where $E_{j}^{l}$ ($l=1,...,h_{j}$) is the direct sum of the
irreducible representations contained in $E_{j}$ that are equivalent to the
same $\left[ \pi _{j,l}\right] \in \widehat{G}$ (see \cite[Section 2.6]%
{Serre}); then, for every $\pi \in \widehat{G}$, $f^{\pi }$ (as defined in (%
\ref{fpi})) is equal to zero if $\left[ \pi \right] \neq \left[ \pi _{j,l}%
\right] $ for every $l=1,...,h_{j}$, and equal to the projection of $f$ on $%
E_{j}^{l}$ if $\left[ \pi \right] =\left[ \pi _{j,l}\right] $ for some $%
l=1,...,h_{j}$.
\end{enumerate}
\end{proposition}

\begin{proof}
\textbf{(1.) }Just write
\begin{eqnarray*}
\int_{T}\left( \mathbb{X}_{N}^{\pi }\left( t\right) -\mathbb{X}^{\pi }\left(
t\right) \right) ^{2}m\left( dt\right) &=&d_{\pi }^{2}\int_{T}\left(
\int_{G}\chi _{\pi }\left( g\right) \left( \mathbb{X}_{N}\left( g^{-1}\cdot
t\right) -\mathbb{X}\left( g^{-1}\cdot t\right) \right) dg\right)
^{2}m\left( dt\right) \\
&\leq &d_{\pi }^{2}\alpha _{\pi }^{2}\int_{T}\left[ \int_{G}\chi _{\pi
}\left( g\right) \left( \left( \mathbb{X}_{N}\left( g^{-1}\cdot t\right) -%
\mathbb{X}\left( g^{-1}\cdot t\right) \right) \right) ^{2}dg\right] m\left(
dt\right) \\
&=&d_{\pi }^{2}\alpha _{\pi }^{2}\int_{T}\left( \left( \mathbb{X}_{N}\left(
t\right) -\mathbb{X}\left( t\right) \right) \right) ^{2}m\left( dt\right)
\rightarrow 0\text{,}
\end{eqnarray*}%
thanks to Assumption C, as well as the $G$-invariance of $m$.

\textbf{(2.) }A function $f$ is in $E_{j}$ if, and only if,
\begin{equation*}
\lambda_j f\left( t\right) =\int_{T}R\left( t,s\right) f\left(
s\right) m\left( ds\right) \text{.}
\end{equation*}%
Now suppose $\mathbb{X}$ has a $G$-invariant law. Then, $R$ is also $G$%
-invariant, and moreover, for every $g\in G$ and $f\in E_{j}$,%
\begin{eqnarray*}
f\left( g^{-1}\cdot t\right) &=&\int_{T}R\left( g^{-1}\cdot t,s\right)
f\left( s\right) m\left( ds\right) \\
&=&\int_{T}R\left( g^{-1}\cdot t,g^{-1}\cdot s\right) f\left( g^{-1}\cdot
s\right) m\left( ds\right) \text{ \ \ (}G\text{-invariance of }m\text{)} \\
&=&\int_{T}R\left( t,s\right) f\left( g^{-1}\cdot s\right) m\left( ds\right)
\text{ \ \ (}G\text{-invariance of }R\text{),}
\end{eqnarray*}%
and therefore $f\in E_{j}$. This concludes the proof.

\textbf{(3.) }This point is a direct application of Theorem 8 in \cite{Serre}%
.
\end{proof}

\subsection{Watson's identity on  the n-dimensional flat tori}
Watson's identity concerns processes defined on $[0,1]$ and taking
the same values at $t=0$ and $t=1$, in other words on a circle.
Among the various geometrical sets arising as generalizations of
the circle in higher dimensions, we will consider the
$n$-dimensional torus. Recall that an $n-$dimensional lattice is a
set
$$\Gamma:=\{\sum_{i=1}^na_iv_i:a_1,...,a_n\in \mathbb{Z}\}$$
where $v_1,...v_n$ are $n$ independent vectors in $\R^n$. The dual
lattice $\Gamma^*$ is defined to be the set of $v^*\in\R^n$ such
that
$$<v|v^*>\in \mathbb{Z},\ \text{for all } v\in\Gamma.$$The
quotient space $T_\Gamma:=\R^n/\Gamma$ is the  $n-$dimensional
torus associated to $\Gamma$, and it is endowed with the measure
$dm$ inherited from the Lebesgue measure on $\mathbb{R}^n$.
Consider a centered Gaussian process $X:=\{X(t):t\in T_\Gamma\}$,
with covariance function $K$. For $n=1$, $\Gamma=\mathbb{Z}$, one
has $T^1=\R/\mathbb Z$ and  $X$ is a centered Gaussian process
defined on $[0,1]$ such that $X(0)=X(1)$. In this case $X$ can be
a Brownian bridge or the Watson process. These processes are
involved in Watson's identity ($\ref{Watson1961}$). We  propose an
assumption on $X$, implying that this process satisfies an
identity analogue of Watson's duplication identity
$(\ref{Watson1961})$ in higher dimensions. The techniques we adopt
represent a $n$-dimensional generalization of the line of reasoning that the second author used in \cite{Pycke}.\\

\bigskip

\textbf{Assumption D} -- There exists a function $k:T_\Gamma\to
\R$ such that
\begin{equation}\label{cov}
K(s,t)=k(t-s)\quad (s,t\in T_\Gamma)
\end{equation}

\bigskip

Note that this assumption is equivalent to the hypothesis that $K$
is invariant under the isometry group of $T_\Gamma$ which is
composed of all translations of vector
$v\in\{\sum_{i=1}^na_iv_i:0\leq a_i\leq 1,1\leq i\leq n\}$. Let us
first check that the covariance function of Watson's process given
by $(\ref{CovCompBr})$ can be put in the form $(\ref{cov})$. If
for $s\in\R$ we denote by $\overline s\in[0,1)$ the corresponding
class in $\R/\mathbb{Z}$, we have, for $s,t\in[0,1]$,
$$|s-t|-\frac12=\begin{cases}t-s-\frac12=\frac12-(1+s-t)=\frac12-\overline{s-t}&\text{if $s< t$,}
\\
s-t-\frac12=\overline{(s-t)}-\frac12&\text{if $s\geq
t$,}\end{cases}$$hence $$|s-t|-\frac12=\pm
\left(\overline{(s-t)}-\frac12\right)\quad (s,t\in[0,1]).$$ This
allows us to obtain the expression \begin{equation}\label{}
s\wedge t-\frac{s+t}{2}+\frac{\left(
s-t\right)^{2}}2+\frac1{12}=\frac12\left(|s-t|-\frac12\right)^2-\frac1{24}=\frac12\left(\overline{s-t}-\frac12\right)^2-\frac1{24}=:k(u)
\end{equation}
where $k(u)=(\overline u-1/2)^2/2-1/24$ for $u\in T^1$.
\begin{lemma}If the centered Gaussian process $X$ satisfies
Assumption $D$, then $K$ admits a Karhunen-Lo\`eve expansion of
the form
\begin{equation}\label{}
    K(s,t)=\sum_{v\in\Gamma^*}\lambda_v\{\alpha_v\cos(2\pi<v|s>)\cos(2\pi<v|t>)+
    \alpha_v\sin(2\pi<v|s>)\sin(2\pi<v|t>)\}
\end{equation}
where $\lambda_v\in[0,\infty)$ for each $v\in\Gamma^*$, and
$\alpha_v>0$ is chosen such that
$$\int_{T_\Gamma}\alpha_v^2\cos^2(2\pi<v|s>)dm(s)=\int_{T_\Gamma}\alpha_v^2\sin^2(2\pi<v|s>)dm(s)=1.$$
\end{lemma}
\begin{proof}The functions $\{u\mapsto \cos(2\pi<v|u>),u\mapsto
\sin(2\pi<v|u>):v\in\Gamma^*\}$ form a complete  set of orthogonal
functions in $L^2(T_\Gamma)$. The Fourier series of $k$ in this
basis has the form
$$k(u)=\sum_{v\in\Gamma^*}\{a_v\cos(2\pi<v|u>)+b_v\sin(2\pi<v|u>)\}.$$
Since $k(u)=k(x-y)=K(x,y)=K(y,x)=k(y-x)=k(-u)$, one has $b_v=0$
for each $v\in\Gamma^*$. If we replace $u$ by $x-y$ and use the
identity $\cos(a-b)=\cos a\cos b+\sin a\sin b$ we obtain the
desired K-L expansion of $K$.
\end{proof}
\begin{theorem} If the centered Gaussian process $X$ satisfies
the assumption $D$, then
\begin{equation}\label{}
    \int_TX^2(t)dm(t)=\frac14\int_TX_1^2(t)dm(t)+\frac14\int_TX_2^2(t)dm(t)
\end{equation}
where $X_1(t):=\frac{X(t)-X(-t)}2$ and
$X_2(t):=\frac{X(t)+X(-t)}2$ are two independent centered Gaussian
processes such that
$$\int_TX_1^2(t)dm(t)\egl\int_TX_2^2(t)dm(t)$$
and $X_1(0,...,0)=X_1(\frac{v_1}2,...,\frac{v_n}2)=0$.
\end{theorem}
\begin{proof}From the preceding Lemma, $X$ has a K-L expansion of the form
$$X(t)=\sum_{v\in\Gamma^*}\lambda_v\{\xi_v\alpha_v\cos(2\pi<v|t>)+\xi'_v\alpha_v\sin(2\pi<v|t>)\}$$
(where the $\xi_v$ and the $\xi_v'$ are independent standard
Gaussian random variables) and the claimed identity is clearly
fulfilled with
$$X_1(t)=\sum_{v\in\Gamma^*}\lambda_v\xi_v\alpha_v\sin(2\pi<v|t>),\
X_2(t)=\sum_{v\in\Gamma^*}\lambda_v\xi_v\alpha_v\sin(2\pi<v|t>)$$
 \end{proof}


\bigskip

\bigskip

\begin{thebibliography}{99}
\bibitem{Adler} Adler R. J. (1990). \textit{An Introduction to Continuity,
Extrema, and Related Topics for General Gaussian Processes}. Lecture
Notes-Monograph Series, \textbf{12}, Institut of Mathematical Statistics,
Hayward, California.

\bibitem{Chou} Chou C.-S. (1994) \textquotedblleft Sur l'extension d'une
identit\'{e} en loi entre le pont brownien et la variance du mouvement
brownien\textquotedblright , \textit{Stochastic processes and their
applications }\textbf{49}, 217-225

\bibitem{DPY} Deheuvels P., Peccati G. and Yor M. (2004), \textquotedblleft
Quadratic functionals of the Brownian sheet and related
processes\textquotedblright , to appear in: \textit{Stochastic Processes and
Their Applications}

\bibitem{DMM} Dellacherie C., Maisonneuve B. and Meyer P.A. (1992). \textit{%
Probabilit\'{e}s et Potentiel, Chapitres XVII \`{a} XXIV}. Hermann, Paris.

\bibitem{Diaconis} Diaconis P. (1988) \textit{Group Representations in
Probability and Statistics. }IMS Lecture Notes -- Monograph Series \textbf{11%
}. Hayward, California

\bibitem{Dudley} Dudley R. (2001) \textit{Real analysis and probability }(%
\textit{2nd edition})\textit{. }Wadsworth and Brooks/Cole, Pacific Grove, CA.%
\textit{\ }

\bibitem{DuiKolk} Duistermaat J.J. and Kolk J.A.C. (1997) \textit{Lie groups}%
. Springer-Verlag. Berlin-Heidelberg-New York

\bibitem{Fox Taqqu} Fox R. and Taqqu M. (1987) \textquotedblleft Multiple
Stochastic Integrals with Dependent Integrators\textquotedblright , \textit{%
Journal of Multivariate Anaysis }\textbf{21}, 105-127

\bibitem{PeYor2005} Peccati G. and Yor M.\ (2005) \textquotedblleft
Identities in law between quadratic functionals of bivariate Gaussian
processes, through Fubini theorems and symmetric
projections\textquotedblright . To appear in the volume: \textit{%
Approximation and Probability}, Banach Center Publications, Warsaw, Poland

\bibitem{Pycke} Pycke J.-R. (2005) \textquotedblleft Sur une identit\'{e} en
loi entre deux fonctionnelles quadratiques du pont
Brownien\textquotedblright . \textit{C. R. Math. Acad. Sci. Paris }(\textit{S%
\'{e}rie I}) \textbf{340} (5), 373-376

\bibitem{PyckeBis} Pycke J.-R. (2005) \textquotedblleft Asymptotic
decompositions for an invariant test of uniformity on the sphere via group
representations\textquotedblright , submitted preprint

\bibitem{RotaShen} Rota G.-C. and Shen J. (2000) \textquotedblleft On the
combinatorics of cumulants\textquotedblright , \textit{Journal of
Combinatorial Theory }(\textit{Series A}), \textbf{91}, 283-304

\bibitem{Rota Wall} Rota G.-C. and Wallstrom T.C. (1997) \textquotedblleft
Stochastic integrals: a combinatorial approach\textquotedblright . \textit{%
Ann. Probab.} \textbf{25}(3), 1257--1283

\bibitem{Serre} Serre J.-P. (1977) \textit{Linear representations of finite
groups. }Graduate Texts in Mathematics \textbf{42}. Springer-Verlag.
Berlin-Heidelberg-New York

\bibitem{ShiYor} Shi Z. and Yor M. (1997) \textquotedblleft On an identity
in law for the variance of the Brownian bridge\textquotedblright . \textit{%
Bull. London Math. Soc. }\textbf{29} (1), 103-108

\bibitem{SW} Shorack, G. R. and Wellner, J. A. (1986) \textit{Empirical
processes with applications to statistics.} Wiley, New York.

\bibitem{Slud} Slud E. V. (1993). \textquotedblleft The moment problem for
polynomial forms in normal random variables\textquotedblright . \textit{Ann.
Probab.} \textbf{21}(4), 2200--2214

\bibitem{Surg} Surgailis, D. (2000) \textquotedblleft CLTs for polynomials
of linear sequences: Diagram formula with illustrations\textquotedblright ,
in: \textit{Long Range Dependence}, 111-128. Birkh\"{a}user, Basel.

\bibitem{Wat} Watson G.S. (1961) \textquotedblleft Goodness-of-fit tests on
a circle\textquotedblright . \textit{Biometrika }\textbf{48}, 109-114
\end{thebibliography}
\end{document}